\documentclass[]{amsart}
\usepackage{amssymb,amsmath}
\usepackage[mathscr]{euscript}
\newcounter{sec}

\newcounter{punct}[sec]

\def\punct{\refstepcounter{punct}{\arabic{sec}.\arabic{punct}.  }}

\newtheorem{theorem}{Theorem}[sec]
\newtheorem{proposition}[theorem]{Proposition}

\newtheorem{lemma}[theorem]{Lemma}

\newtheorem{corollary}[theorem]{Corollary}

\def\COUNTERS{\addtocounter{sec}{1}
              \setcounter{punct}{0}
          \setcounter{equation}{0}
          \setcounter{theorem}{0}
         }
          
          \def\sm{\smallskip}
          
\begin{document}

\newcommand{\supp}{\mathop {\mathrm {supp}}\nolimits}
\newcommand{\rk}{\mathop {\mathrm {rk}}\nolimits}
\newcommand{\Aut}{\mathop {\mathrm {Aut}}\nolimits}
\newcommand{\Out}{\mathop {\mathrm {Out}}\nolimits}
\newcommand{\OO}{\mathop {\mathrm {O}}\nolimits}
\renewcommand{\Re}{\mathop {\mathrm {Re}}\nolimits}
\newcommand{\ch}{\cosh}
\newcommand{\sh}{\sinh}
\newcommand{\Op}{\mathop {\mathrm {Op}}\nolimits}
\newcommand{\tr}{\mathop {\mathrm {tr}}\nolimits}
\newcommand{\diag}{\mathop {\mathrm {diag}}\nolimits}

\def\0{\mathbf 0}

\def\ov{\overline}
\def\wh{\widehat}
\def\wt{\widetilde}

\renewcommand{\rk}{\mathop {\mathrm {rk}}\nolimits}
\renewcommand{\Aut}{\mathop {\mathrm {Aut}}\nolimits}
\renewcommand{\Re}{\mathop {\mathrm {Re}}\nolimits}
\renewcommand{\Im}{\mathop {\mathrm {Im}}\nolimits}
\newcommand{\sgn}{\mathop {\mathrm {sgn}}\nolimits}
\newcommand{\Isoc}{\mathop {\mathrm {Isoc}}\nolimits}
\newcommand{\PIsoc}{\mathop {\mathrm {PIsoc}}\nolimits}

\newcommand{\Sch}{\mathop {\mathrm {Sch}}\nolimits}
\newcommand{\sch}{\mathop {\mathrm {sch}}\nolimits}
\newcommand{\Fr}{\mathop {\mathrm {Fr}}\nolimits}

\def\bfa{\mathbf a}
\def\bfb{\mathbf b}
\def\bfc{\mathbf c}
\def\bfd{\mathbf d}
\def\bfe{\mathbf e}
\def\bff{\mathbf f}
\def\bfg{\mathbf g}
\def\bfh{\mathbf h}
\def\bfi{\mathbf i}
\def\bfj{\mathbf j}
\def\bfk{\mathbf k}
\def\bfl{\mathbf l}
\def\bfm{\mathbf m}
\def\bfn{\mathbf n}
\def\bfo{\mathbf o}
\def\bfp{\mathbf p}
\def\bfq{\mathbf q}
\def\bfr{\mathbf r}
\def\bfs{\mathbf s}
\def\bft{\mathbf t}
\def\bfu{\mathbf u}
\def\bfv{\mathbf v}
\def\bfw{\mathbf w}
\def\bfx{\mathbf x}
\def\bfy{\mathbf y}
\def\bfz{\mathbf z}

\def\bfA{\mathbf A}
\def\bfB{\mathbf B}
\def\bfC{\mathbf C}
\def\bfD{\mathbf D}
\def\bfE{\mathbf E}
\def\bfF{\mathbf F}
\def\bfG{\mathbf G}
\def\bfH{\mathbf H}
\def\bfI{\mathbf I}
\def\bfJ{\mathbf J}
\def\bfK{\mathbf K}
\def\bfL{\mathbf L}
\def\bfM{\mathbf M}
\def\bfN{\mathbf N}
\def\bfO{\mathbf O}
\def\bfP{\mathbf P}
\def\bfQ{\mathbf Q}
\def\bfR{\mathbf R}
\def\bfS{\mathbf S}
\def\bfT{\mathbf T}
\def\bfU{\mathbf U}
\def\bfV{\mathbf V}
\def\bfW{\mathbf W}
\def\bfX{\mathbf X}
\def\bfY{\mathbf Y}
\def\bfZ{\mathbf Z}
\def\bfT{\mathbf T}

\def\frD{\mathfrak D}
\def\frX{\mathfrak X}
\def\frS{\mathfrak S}
\def\frZ{\mathfrak Z}
\def\frL{\mathfrak L}
\def\frG{\mathfrak G}
\def\frg{\mathfrak g}
\def\frh{\mathfrak h}
\def\frf{\mathfrak f}
\def\frl{\mathfrak l}
\def\frp{\mathfrak p}
\def\frq{\mathfrak q}
\def\frr{\mathfrak r}

\def\bfw{\mathbf w}
%%% END MATHBF
%%%%%%%%%%%%%%%%%%%%%%%%%%%%%%%
%%%%%%%%%%%%%%%%%%%%%%%%%%%%%%%%%
%%% BEGIN MATHBB

\def\R {{\mathbb R }}
 \def\C {{\mathbb C }}
  \def\Z{{\mathbb Z}}
\def\K{{\mathbb K}}
\def\N{{\mathbb N}}
\def\Q{{\mathbb Q}}
\def\A{{\mathbb A}}
\def\U{{\mathbb U}}

\def\T{\mathbb T}
\def\P{\mathbb P}

\def\G{\mathbb G}

\def\cD{\EuScript D}
\def\cL{\mathscr L}
\def\cK{\EuScript K}
\def\cM{\EuScript M}
\def\cN{\EuScript N}
\def\cP{\EuScript P}
\def\cQ{\EuScript Q}
\def\cR{\EuScript R}
\def\cT{\EuScript T}
\def\cW{\EuScript W}
\def\cY{\EuScript Y}
\def\cF{\EuScript F}
\def\cG{\EuScript G}
\def\cZ{\EuScript Z}
\def\cI{\EuScript I}
\def\cB{\EuScript B}
\def\cA{\EuScript A}
\def\cO{\EuScript O}
\def\cE{\EuScript E}

\def\bbA{\mathbb A}
\def\bbB{\mathbb B}
\def\bbD{\mathbb D}
\def\bbE{\mathbb E}
\def\bbF{\mathbb F}
\def\bbG{\mathbb G}
\def\bbI{\mathbb I}
\def\bbJ{\mathbb J}
\def\bbL{\mathbb L}
\def\bbM{\mathbb M}
\def\bbN{\mathbb N}
\def\bbO{\mathbb O}
\def\bbP{\mathbb P}
\def\bbQ{\mathbb Q}
\def\bbS{\mathbb S}
\def\bbT{\mathbb T}
\def\bbU{\mathbb U}
\def\bbV{\mathbb V}
\def\bbW{\mathbb W}
\def\bbX{\mathbb X}
\def\bbY{\mathbb Y}

\def\kappa{\varkappa}
\def\epsilon{\varepsilon}
\def\phi{\varphi}
\def\le{\leqslant}
\def\ge{\geqslant}

\def\B{\mathrm B}

\def\la{\langle}
\def\ra{\rangle}
\def\tri{\triangleright}

\def\lambdA{{\boldsymbol{\lambda}}}
\def\alphA{{\boldsymbol{\alpha}}}
\def\betA{{\boldsymbol{\beta}}}
\def\mU{{\boldsymbol{\mu}}}

\def\const{\mathrm{const}}
\def\rem{\mathrm{rem}}
\def\even{\mathrm{even}}
\def\SO{\mathrm{SO}}
\def\SL{\mathrm{SL}}
\def\SU{\mathrm{SU}}
\def\GL{\operatorname{GL}}
\def\End{\operatorname{End}}
\def\Mor{\operatorname{Mor}}
\def\Aut{\operatorname{Aut}}
\def\inv{\operatorname{inv}}
\def\red{\operatorname{red}}
\def\Ind{\operatorname{Ind}}
\def\dom{\operatorname{dom}}
\def\im{\operatorname{im}}
\def\md{\operatorname{mod\,}}
\def\St{\operatorname{St}}
\def\Ob{\operatorname{Ob}}
\def\PB{{\operatorname{PB}}}
\def\Tra{\operatorname{Tra}}

\def\ZZ{\mathbb{Z}_{p^\mu}}
\def\F{\mathbb{F}}

\def\cH{\EuScript{H}}
\def\cQ{\EuScript{Q}}
\def\cL{\EuScript{L}}
\def\cX{\EuScript{X}}

\def\Di{\Diamond}
\def\di{\diamond}

\def\fin{\mathrm{fin}}
\def\ThetA{\boldsymbol {\Theta}}

\def\0{\boldsymbol{0}}

\def\FF{\,{\vphantom{F}}_3F_2}
\def\HH{\,\vphantom{H}^{\phantom{\star}}_3 H_3^\star}
\def\Ho{\,\vphantom{H}_2 H_2}

\def\disc{\mathrm{disc}}
\def\cont{\mathrm{cont}}

\def\fan{\vphantom{|^|}}

\def\osigma{\ov\sigma}
\def\ot{\ov t}

\begin{center}
\Large\bf

Oligomorphic groups, categories of partial bijections,
 and ultrahomogeneous cubic spaces over  finite fields

%On the group of symmetries of ultrahomogeneous linear space over a finite field
%equipped with a cubic form

\medskip

\sc

Yury A. Neretin%
\footnote{The work is supported by the grants of FWF (the Austrian Scientific Funds),  P31591, PAT5335224.}

\end{center}

\bigskip

{\small

We show that for a certain class of oligomorphic groups there is a  version of multiplication  of double cosets in the  Ismagilov--Olshanski sense. Categories of (reduced) double cosets are realized as certain categories of partial bijections. 
As an example, we consider the ultrahomogeneous linear space $\EuScript{Z}$ of countable dimension  over a
finite field  equipped with a  cubic form.
The group  of automorphisms of $\EuScript{Z}$ is an oligomorphic group, we 
describe its open subgroups.
According Tsankov, this gives
a classification of its unitary representations. The category of reduced double cosets
in this case is  the category of partial linear bijections of finite-dimensional cubic spaces  preserving cubic forms.

}

\bigskip

\setcounter{sec}{-1}

\section*{0\phantom{a} Introduction}

\COUNTERS

{\bf\punct Categories of double cosets.}
In theory of unitary representations of 'large' ('infinite-dimensional' groups%
\footnote{A.~M.~Vershik proposed the term '{\it large groups}'
	instead of 'infinite-dimensional groups'. Numerous 'infinite-dimensional
	groups' are totally disconnected, hence they have dimension 0 as topological
	spaces. This is so for groups of types $2^\circ$ and $4^\circ$-$7^\circ$ below. We emphasize
	that  infinite dimensional unitary groups and  (zero-dimensional) infinite symmetric groups are close relatives (at least from the point of view
	of representation theory).}), quite often
the following phenomena
 arise. Let $G$ be a group, and $K$ is a subgroup.
 It can happen that the  space of double cosets $K\backslash G/K$
 admits a natural associative multiplication
  $$
  K\backslash G/K\, \times\, K\backslash G/K\, \to \, K\backslash G/K,
  $$
and the semigroup $K\backslash G/K$ acts in  unitary representations
of $G$. 

Initially, such semigroups  appeared as tools for classification of unitary representations
and descriptions of spherical functions.
The first examples  were discovered by Rais Ismagilov \cite{Ism},
\cite{Ism+}. A wider classes of such multiplications appeared in works by
Grigori Olshanski (see \cite{Olsh-sym}--\cite{Olsh-semigroup}).
Now it is clear that this phenomenon is widespread. In particular,
it  arises for the following types of groups:

\sm 

$1^\circ.$ infinite-dimensional real classical groups, see, for example,
\cite{Olsh-GB}, \cite{Olsh-semigroup}, \cite{Ner-book}, \cite{Ner-coll};

\sm 

$2^\circ.$ various versions of 'infinite symmetric groups', see, for example,
\cite{Olsh-sym}, \cite{Olsh-bisym}, \cite{Olsh-semigroup}, 
\cite{Ner-book}, \cite{Ner-field};

\sm

$3^\circ.$ groups of transformations of measure spaces, see, for example, \cite{Ner-book}, \cite{Ner-sin};

\sm 

$4^\circ.$ infinite-dimensional classical groups over finite fields, see, for example, \cite{Olsh-semigroup}, \cite{Ner-finite};

\sm

$5^\circ.$ infinite-dimensional groups over non-Archimedean locally compact fields, see, for example,   \cite{Ner-p-adic};

\sm

$6^\circ.$  groups of automorphisms and spheromorphisms of non-locally finite trees and $\R$-trees, see, for example,
\cite{Olsh-new}, \cite{Ner-urysohn}.

\sm 

{\sc Remarks.}  %a) The references are related to purposes of the present paper and are not well-representative,  see more  in  \cite{Ner-sin}.
a) In all cases, $K$ looks as an infinite-dimensional counterpart of compact groups. In 
\cite{Ner-book}, Chapter VIII, there was proposed a heuristic notion of 
a '{\it heavy group}'.  It turns out  that 
for groups $G$ with a heavy subgroup $K$ the multiplication of double cosets arises
more or less automatically%
\footnote{Examples of heavy groups are the complete symmetric group
$S(\infty)$, the complete infinite hyperoctahedral group, the complete
infinite-dimensional  unitary group $\mathrm{U}(\infty)$, the complete orthogonal group $\mathrm{O}(\infty)$, the group of measure preserving transformations
of a Lebesgue measure space. Oligomorphic groups discussed below look
like 'heavy groups'.
\newline
All such groups $K$ are Roelcke precompact. This means that
for any open subset $U$ in $K$ there is a finite number of elements
$h_j$ such that the sets $U h_j U$ cover $K$ (see, e.g., \cite{Tsa},
\cite{BIT}). All 'heavy groups',
that appeared in such representation-theoretical considerations 
 are Roelcke precompact. But according \cite{Ben},
 the group of all non-singular transformations of a Lebesgue measure space
 is Roelcke precompact, and it looks completely differently.}.

%\sm

%b) As far as I know, for all existing substantial theories of unitary representation
%of real classical groups and infinite symmetric groups, the multiplication of
%double cosets is present (but this is not so for groups of matrices over a finite field; namely,
%this is not so for the construction discussed in  \cite{GKV});

\sm 

b) For groups of matrices over finite fields and local fields existing picture is relatively poor (comparing $1^\circ$--$2^\circ$).
But it is non-trivial and drastically differs from the real case.
\hfill $\boxtimes$

\sm 

It turns out that in all examined cases, we actually have families of subgroups $K_\alpha$
enumerated by elements $\alpha$ of some set $A$. We have
a collection of operations
$$
K_\alpha\backslash K/K_\beta\,\times K_\beta \backslash K/ K_\gamma\, \to\, K_\alpha \backslash K/ K_\gamma.
$$
In fact, we come to a category ({\it train of $G$}), whose objects are enumerated by elements $\alpha\in A$,  
sets of morphisms $\beta\to\alpha$ are  $K_\alpha\backslash K/K_\beta$, and products
of morphisms are products of double cosets.

\sm 

{\bf\punct Oligomorphic groups.} 
The first purpose of the present paper is to add
the line:

\sm 

$7^\circ.$  oligomorphic groups

\sm 

\noindent
to the list $1^\circ$--$6^\circ$  above.
\par

Denote by $S(\Omega)$ the group of all permutations of a countable space $\Omega$.
An {\it oligomorphic group} (see \cite{Cam}) is a closed subgroup in $S(\Omega)$,
which has a finite number of orbits on any finite  product $\Omega\times\dots\times \Omega$. See a collection of
 examples below in Subsect. 
\ref{ss:examples}.
%of finite numbers of copies of $\Omega$.

This notion arose in model theory (which is a branch of mathematical logic). It is interesting to link this topic with other 
investigations in large groups. Some oligomorphic groups
look exotic from the point of view of representation theory, but some of them are important (or at least look  naturally).
 Todor Tsankov \cite{Tsa} proved that all
 unitary representations of oligomorphic groups are induced from finite-dimensional
representations of open subgroups.

This statement allows us to easily describe multiplications of {\it reduced} double cosets for oligomorphic groups (under some additional conditions).
The corresponding categories are certain categories of partial bijections
(the term 'reduced double cosets' is explained below in Subsect. \ref{ss:reduced}). 

The work \cite{Tsa} was used earlier  in \cite{Ner-finite}, \cite{Ner-p-adic}.
Since oligomorphic groups look like heavy groups,
the next (open) question is a search for interesting
overgroups of oligomorphic groups, see a discussion in Subsect. \ref{ss:question}.

%We recall  two examples of oligomorphic groups: the group of linear operators in the countable-dimensional
%linear space over a finite field and the group of automorphisms $\Aut(\cR)$ of the universal Rado graph (for its definition, see below Subsect. \ref{ss:examples}).

% The group $\Aut(\cR)$
%is exotic from the point of view of representation theory. However, it is interesting that such groups
%can be involved to representation theory. For instance,
%the  category of reduced double cosets corresponding to $\Aut(\cR)$ is described 
%in the following way. Its objects are finite graphs. A morphism $X\to Y$ is an isomorphism
%of an induced subgraph%
%\footnote{To avoid a misleading we use the term '{\it induced subgraph}': we consider some set of vertices of
%	$X$ and {\it all} edges of $X$ connecting these vertices.} $A\subset X$ to an induced subgraph $B\subset Y$.
%$*$-Representations of this category (see formal definitions in Subsect. \ref{ss:rep-cat}) are in a canonical one-to-one correspondence with unitary representations
%of the group $\cR$. 

\sm

{\bf\punct Universal cubic spaces and groups of their automorphisms.}
The second purpose of the paper is to obtain a classification of unitary representations
of the group of automorphisms of a universal cubic space over a finite field.

Consider a finite or countable field $\Bbbk$. We say that a {\it cubic space}
is a linear space $Z$ over $\Bbbk$ equipped with a symmetric trilinear
 form
 $Z\times Z\times Z\to \Bbbk$.
A linear space $\cZ(\Bbbk)$ over $\Bbbk$ of countable dimension is
{\it universal} if any finite-dimensional cubic space $X$ over $\Bbbk$ 
admits an embedding to $\cZ(\Bbbk)$, and the set of such embeddings is a homogeneous  space
with respect to the group $\Isoc\bigl(\cZ(\Bbbk) \bigr)$ 
of automorphisms of $\cZ(\Bbbk)$. 
Such space exists%
\footnote{Universal cubic spaces (and more generally universal tensor spaces \cite{HS})
	are representatives of numerous counterparts of the Urysohn universal metric
	spaces \cite{Ury} (as the Rado graph, the Hall universal group, the universal poset, the Gurarii space, ultrametric Urysohn spaces, etc.), see, e.g., \cite{Macph}.} and is unique. It  was a topic of numerous investigations 
in last years, see, for example, \cite{Bik1}--\cite{Bik2}, \cite{KZ},  \cite{HS}.

If a field $\Bbbk$ is finite, then the group $\Isoc\bigl(\cZ(\Bbbk) \bigr)$
is oligomorphic. We classify its open subgroups and unitary representations. The objects of the corresponding category 
are finite dimensional cubic spaces over $\Bbbk$, a morphism $V\to W$
is a form preserving invertible linear operator from a  subspace $A\subset X$ to a subspace $B\subset Y$.

\section{Multiplications of double cosets and trains of 'large' groups}

\COUNTERS

{\bf\punct Multiplication of double cosets.}
%{\bf \punct Multiplication of double cosets and $(G,K)$-pair.}
Let $G$ be a group, $K$,  $K_1$, $K_2$ be  its subgroups. Denote by 
$K_1\backslash G/K_2$ the set of double cosets of $G$ with respect to $K_1$ and $K_2$, i.e., the quotient of $G$ by the equivalence relation
$$
g\sim h_1 g h_2, \quad\text{where $h_1\in K_1$, $h_2\in K_2$.}
$$
We denote classes of equivalence (points of the quotient) by $K_1 g K_2$.

\sm

Let $\rho$ be a unitary representation of a {\it separable} topological
 group $G$ in a {\it separable} Hilbert space $\cL$. Denote
by $\cL^K\subset \cL$ the subspace of $K$-fixed vectors. By $P^K$ we denote the orthogonal projection to $\cL^K$. It is easy to see that an operator
$$
\wh \rho (g) := P^K \rho (g)\, P^K
$$ 
depends only on the double coset $\frg=KgK$ containing $g$. Indeed, let  $k_1$, $k_2\in K$, $\xi$, $\eta\in \cL$.
Then matrix elements of operators $P^K \rho (k_1 g k_2) P^K$ and $P^K \rho (g) P^K$ coincide. Indeed,
\begin{multline*}
\bigl\la P^K \rho (k_1 g k_2) P^K \xi,\,\eta\bigr\ra_\cL= 
\bigl\la P^K \rho (k_1) \rho (g) \rho (k_2) P^K \xi,\,\eta\bigr\ra_\cL=\\=
\bigl\la P^K \rho (k_1) \rho (g)  P^K \xi,\,\eta\bigr\ra_\cL
=\bigl \la \rho (k_1) \rho (g)  P^K \xi,\, P^K  \eta\bigr\ra_\cL=\\=\bigl \la  \rho (g)  P^K \xi,\,\rho (k_1^{-1})P^K  \eta\bigr\ra_\cL=
\bigl \la  \rho (g)  P^K \xi,\,P^K  \eta\bigr\ra_\cL= \bigl \la P^K  \rho (g)  P^K \xi, \,\eta\bigr\ra_\cL.
\end{multline*}

\sm

{\sc Remark.}
Decompose $\cL=\cL^K\oplus (\cL^K)^\bot$. Then the operator 
$\wh \rho (g)$ has a block form 
$\begin{pmatrix}\wt \rho(g)&0\\0&0\end{pmatrix}$, where
$\wt\rho(g)$ is an operator in $\cL^K$.
\hfill $\boxtimes$

\sm

 If a group $G$ is finite,
then for certain non-negative rational constants 
$c_{\vphantom{|}\frp\,\frq}^{\vphantom{|}\frr}$ we have 
\begin{equation}
\wh\rho(\frp)\,\wh\rho(\frq)=\sum_{\frg\in K\setminus G/K} 
c_{\vphantom{|}\frp\,\frq}^{\vphantom{|}\frr}\,\, \wh \rho(\frr)
\label{eq:convolution}
\end{equation}
for any unitary  representation $\rho$ of $G$.

Equivalently, let us consider the group algebra $\C[G]$ of $G$,
denote by $\phi \star \psi$ the convolution in $\C[G]$.
For $p\in G$ denote by $\delta_p\in \C[G]$ the  function,
which is 1 at $p$ and zero at all other points. So,
$\delta_p\star \delta_q=\delta_{pq}$.

Consider the subalgebra $\C[K\setminus G/K]\subset \C[G]$
consisting of $K$-biinvariant functions $\phi$: 
$$ \phi(k_1 g k_2)=\phi(g)\qquad \text{for any $g\in G$,
and $k_1$, $k_2\in K$.} 
$$
For each
 double coset $\frp$, we consider the 'delta-function' of $\frp$:
$$
\Delta_\frp (g)=
\frac 1{\# \frp}\sum_{g\in \frp} \delta_g=
\begin{cases}
\frac 1{\# \frp}, \quad &\text {if $g\in \frp$};
\\
0,\quad &\text{otherwise},
\end{cases}
$$
where $\# \frp$ denotes the number of elements of $\frp$.
Such functions form a basis in $\C[K\setminus G/K]$.
The function $\Delta_K:=\Delta_{KeK}$ is an idempotent
in $\C[G]$ and the unit in $\C[K\setminus G/K]$:
$$
\Delta_K\star \Delta_K=\Delta_K,
\quad \Delta_K \star \delta_p \star \Delta_K=\Delta_{KpK},
\quad \Delta_K \star \Delta_\frp=\Delta_\frp= \Delta_\frp \star \Delta_K. 
$$
For any unitary representation $\rho$ of $G$
we have 
$$
\rho\bigl(\Delta_K\bigr)=P^K, 
$$
 and therefore, 
 $$
 \rho \bigl(\Delta_{KpK}\bigr)=\rho\bigl(\Delta_K \star \delta_p \star \Delta_K\bigr)
= \rho\bigl(\Delta_K \bigr)\, \rho\bigl(\delta_p) \, \rho\bigl(\Delta_K \bigr) 
=P^K \, \rho\bigl(p) \, P^K
=\wh\rho(\frp).
 $$
 So, the coefficients in \eqref{eq:convolution}
 are simply the structure constants of the convolution algebra $\C[K\backslash G/K]$,
 \begin{equation}
 \Delta_\frp\star \Delta_\frq=\sum_{\frg\in K\setminus G/K} 
c_{\vphantom{|}\frp\,\frq}^{\vphantom{|}\frr}\,\, \Delta_\frr.
\label{eq:structure}
\end{equation}

 If $G$ is locally compact and $K$ is a compact subgroup, then there is 
a system of probabilistic measures $\mu_{\fan\frp\,\frq}$ on $K\setminus G/K$ such that 
\begin{equation}
\wh\rho(\frp)\,\wh\rho(\frq)=
\int_{\frr\,\in K\setminus G/K}
 \wh \rho(\frr)\, d\mu_{\fan\frp\,\frq}\,(\frr).
 \label{eq:convolution2} 
\end{equation}
for all unitary representations $\rho$.

Equivalently, we can consider convolution algebras of $K$-biinvariant
functions (or $K$-biinvariant measures) on $G$.

 Such '{\it Hecke--Iwahori algebras}' of double cosets were widely used in classical 
representation theory. But except few examples (including Hecke algebras, affine Hecke 
algebras%
\footnote{Both objects and names for them were introduced by Nagayoshi Iwahori.},
and closely related Yokonuma algebras) algebras of double cosets are heavy
objects.

It turns out that for a 'large group' $G$ and its subgroup $K$
 quite often
 there exists 
a natural associative multiplication $(\frp,\frq)\mapsto \frp\circ\frq$ on
the spae $K\backslash G/K$ such that
 \begin{equation}
 \wh\rho(\frp)\,\wh\rho(\frq)=
 \wh\rho(\frp\circ\frq)
 \label{eq:multiplicativity0}
\end{equation}
for all unitary representations $\rho$ of the group $G$.
% (see, e.g., Grigori Olshanski \cite{Olsh-sym}--\cite{Olsh-semigroup} and the author \cite{Ner-book}--\cite{Ner-sin},
%the first constructions of this type were obtained by Rais Ismagilov \cite{Ism},
%\cite{Ism+}).

\sm

{\sc Remark.}
%In an informal sence, this means that if $p\in \frp$, $q\in \frq$, then %$pq$ almost certainly is contained
%in one double coset. But 'almost certainly' here is completele heuristic. 
Olshanski
\cite{Olsh-sym}, \cite{Olsh-GB}, \cite{Olsh-semigroup} proposed
the following explanation of this phenomenon. Let $G_n$ be finite-dimensional (or finite) groups, let $G_\infty$ be their 'limit'
(as symmetric groups $S_n$ and $S_\infty$).  
Then, in formulas of type \eqref{eq:structure},  we get a concentration of convolutions near one double coset
% see a wider discussion in \cite{Ner-concentration}
(but this explanation does not always works; moreover;   also,
 this does not
formally implies \eqref{eq:multiplicativity0}).
\hfill $\boxtimes$

\sm 

 It is important that such multiplications $\frp\circ\frq$ admit transparent
descriptions (in a contrary with formulas of the type
\eqref{eq:convolution}--\eqref{eq:convolution2}).
%  Often such multiplications lead to unusual algebraic structures, which can be interesting by themselves (see Olshanski
%\cite{Olsh-bisym}--\cite{Olsh-semigroup} and the author \cite{Ner-book}, Sect. IX.3-4, 
%\cite{Ner-coll},
%\cite{Ner-field}). 

  \sm

{\bf \punct Reduced double cosets.%
\label{ss:reduced}}
Generally, functions $\frg\mapsto\wh \rho(\frg)$, where $\rho$ ranges in the space of all unitary
representations, do not separate
 elements  $\frg$ of $K\backslash G/K$. Definitely, this is so if the
 topological quotient space $K\backslash G/K$
 is non-Hausdorff (orbit spaces of non-compact groups are non-Hausdorff very often).
  Olshanski (\cite{Olsh-semigroup}, Sect. 4) found an example,
  when this  happens  
  for finite 
  Hausdorff spaces $K\backslash G/K$.  
  
  \sm
  
  Consider a separable topological  group $G$, its subgroup 
$K$, and a family of subgroups $K_\alpha\subset K$ enumerated by
a set $\cA$ (we assume that this family includes $K$).
For any unitary representation $\rho$ of $G$ in a Hilbert apace $\cL$, we denote 
by $\cL_\alpha$ the subspace of $K_\alpha$-fixed vectors, by $P_\alpha$ the orthogonal projection to $\cL_\alpha$. Consider operators
$$
\wt \rho_{\alpha, \beta}(g):= P_\alpha\,\rho(g)\Bigr|_{\cL_\beta}:
\, \cL_\beta\to \cL_\alpha.
$$

We define {\it reduced double cosets} $[\![K_\alpha\backslash G/K_\beta]\!]$
 as quotients of $G$ by the following equivalence relation:
  $$
  g_1\sim g_2 \quad \text{if $\wt\rho(g_1)=\wt\rho(g_2)$
   for all unitary representations $\rho$.}
  $$
We say that $G$ is a {\it $(G,K)$-pair} if for each $\alpha$, $\beta$, 
$\gamma\in\cA$ for any $\frp\in [\![K_\alpha\backslash G/K_\beta]\!]$,
$\frq\in[\![K_\beta\backslash G/K_\gamma]\!]$ there exists an element
$\frp\circ \frq \in[\![K_\alpha\backslash G/K_\gamma]\!]$ such that
the following property of {\it multiplicativity} holds:
\begin{equation}
\wt\rho_{\alpha,\beta}(\frp)\,\wt\rho_{\beta,\gamma}(\frq)=\wt\rho_{\alpha,\gamma}(\frp\circ \frq )
\label{eq:multiplicativity}
\end{equation}
for all unitary representations $\rho$ of $G$.

Since a product of linear operators is an associative operation,
 the product $\frp\circ \frq $ also is associative. So we come to a category
$\Tra\bigl(G,\{K_\alpha\}\bigr)$ (the {\it train of $G$}). Its objects are enumerated by elements
of $\cA$. Morphisms $\beta\to\alpha$ are elements of 
$[\![K_\alpha\backslash G/K_\beta]\!]$. A product of morphisms is
the $\circ$-product.

This category is equipped with the
{\it involution} $\frp\mapsto \frp^*$ acting from
$[\![K_\alpha\backslash G/K_\beta]\!]$ to $[\![K_\beta\backslash G/K_\alpha]\!]$ for each pair $\alpha$, $\beta$,
it is
determined by the map $g\mapsto g^{-1}$. We have 
$$(\frp\circ \frq)^*=\frq^*\circ \frp^*, \qquad \frp^{**}=\frp.$$

By the definition, any unitary representation $\rho$ of the group $G$ generates a representation $\wt \rho$
 of the category  $\Tra(G,\{K_\alpha\})$ (for a formal definition of  representations
 of categories, see the next subsection).
 
% \sm
 
% {\sc Remark.} Collection of constructions in this spirit is huge.
% Products of double cosets and multiplicativity exist
% for various versions of infinite symmetric groups,
% for classical groups over $\R$, over finite fields, over $p$-adics,
% for groups of transformations of measure spaces, for groups
% acting on non-locally finite trees and $\R$-trees
% (see \cite{Ner-book}--\cite{Ner-sin}, \cite{Olsh-sym}, \cite{Olsh-semigroup}, see more
%formal comments in \cite{Ner-sin}). In this paper, we add oligomorphic groups (see the next section) to this 
%zoo.  \hfill $\boxtimes$
 
 \sm

{\bf\punct Representations of categories.%
\label{ss:rep-cat}}
Recall basic definitions and notation related to representations of categories.
For details, discussion, and further definitions, see \cite{Ner-book}, Sect. II.8.

Let $\cK$ be  a  category. 
Denote by $\Ob(\cK)$ the class of  objects of  $\cK$, by
$\Mor(V,W)=\Mor_{\fan\cK}(V,W)$ the collection of all morphisms $V\to W$, by $\End(V)$ the semigroup
of all endomorphisms of an object $V$. We assume that these semigroups have 
units $1_V$ and these units are units in the category $\cK$, i.e.,
for any $\frp\in \Mor(V,W)$ we have
$$
\frp\cdot 1_V=1_W\cdot \frp =\frp.
$$
Denote by $\Aut(V)$ the group of automorphisms of an object $V$.
An {\it involution} in $\cK$ is a collection of maps $\frp\mapsto\frp^*$
acting from $\Mor(V,W)$ to $\Mor(W,V)$ for all $V$, $W$
 such that
$(\frp\frq)^*=\frq^*\frp^*$, $\frp^{**}=\frp$. 

A 
{\it $*$-representation $(Q,\kappa)$ of a  category} $\cK$ 
is a functor from $\cK $ to the category of separable Hilbert spaces and bounded linear operators.
Namely, for each $V\in\Ob(\cK)$ we assign a Hilbert space $Q(V)$
and for each $\frp\in \Mor(V,W)$ we assign a bounded operator 
$$\kappa_{\fan W,V}(\frp):Q(V)\to Q(W)$$
such that for any $V$, $W$, $Y\in\Ob(\cK)$ and for any $\frp\in\Mor(V,W)$, $\frq\in \Mor(W,Y)$ we have 
$$\kappa_{\fan Y,W}(\frq)\,\kappa_{\fan W,V}(\frp)=\kappa_{\fan Y,V}(\frq\frp).$$
 We also
require 
$$\kappa_{\fan V,W}(\frp)^*=\kappa_{\fan W,V}(\frp^*), \qquad
\kappa_{\fan V,V}(1_V)=1_{\fan Q(V)}
.$$

Consider two $*$-representations $(Q,\kappa)$ and $(Q',\kappa')$ of $\cK$.
An {\it intertwiner} $A:(Q,\kappa)\to (Q',\kappa')$
is a collection
 of 
bounded operators $A(V):Q(V)\to Q'(V)$  
 such that
for any $\frp\in\Mor(V,W)$ we have
$$
\kappa_{W,V}'(\frp) A(V)=A(W) \kappa_{W,V}(\frq).
$$

Two $*$-representations $(Q,\kappa)$ and $(Q',\kappa')$ 
are {\it equivalent} if there is an intertwiner $A$ such
that all operators $A(V)$ are unitary.
 
\sm 
 
 {\sc Remark.}  %We must say some formalities.
We consider representations of categories up to equivalence. For this reason, for us, there is no difference between equivalent categories.
 Recall that according Bernays--G\"odel version of set theory, a collection of objects
of a category is a class (and, generally, not a set).
{\it By definition, all categories $\Tra(G,\{K_\alpha\})$
  are small.} This means that collections of objects and
   collections of morphisms of our categories are sets (and not  `classes').
By aesthetic reason, sometimes we consider  {\it essentially
  small categories}, i.e., categories equivalent to small categories.
 \hfill $\boxtimes$   
  
\sm

{\bf \punct Category of partial bijections.%
\label{ss:PB}}
   Recall the definition of   partial bijections
(see, e.g., \cite{Ner-book}, Sect. VIII.1). Let $X$, $Y$ be sets.
A {\it partial bijection} $\mu: X\to Y$ is a bijection 
from a subset $A\subset X$ to a subset $B\subset Y$. We say that
$A$ is the {\it domain} $\dom \mu$ of $\mu$, $B$ is the {\it image} $\im \mu$.
A {\it product} of partial bijections $\mu:X\to Y$,
 $\nu:Y\to Z$ is defined in the following way:  $\dom\nu\mu:=\mu^{-1}(\dom \nu\,\cap\, \im \mu)$ and $\nu\mu(x)=\nu(\mu(x))$
for $x\in\dom \nu\mu$.

 We also define the {\it pseudoinverse} partial bijection, 
 if $\mu:X\to Y$ is as above,
 then its pseudoinverse  $\mu^*:X\to Y$ is the inverse bijection $B\to A$. 

Clearly, we get a category with involution. Below we consider the category $\PB$
of partial bijections of finite sets.
Clearly, $\PB$ is  essentially small.
For instance, we come to an equivalent small category, if we fix a countable set $\Omega$ and reduce a collection
of objects to the set of all finite subsets of $\Omega$ (morphisms are 
the same).
On the other hand, we can choose a set of each finite
cardinality, and again we come to a small category.

%Denote by $\PB$ the category of partial bijections of finite sets.

\sm
 
 {\bf \punct Example: the infinite symmetric group $S(\infty)$
 and the category of partial bijections.%
 \label{ss:S-infty}} 
 See Olshanski \cite{Olsh-sym}, see also \cite{Ner-book}, Sect. VIII.1-2.
 Let $\Omega$ be a countable or finite  set, $S(\Omega)$ be the group of {\it all} permutations 
 of $\Omega$. Denote $S_\infty:=S(\N)$.
  It is known, that this group ({\it infinite symmetric group}) has a unique separable topology (see \cite{KR}). Namely, $g_j$ converges to $g$ if for each $\omega\in \Omega$ we have $g_j\omega=g\omega$ for sufficiently large $j$.
  
  Let $G=K=S_\infty$. Denote by $K_{\alpha}\subset S_\infty$
  the subgroup consisting of elements fixing $1$, \dots, $\alpha\in \N$,
  we assume
  $K_0:=S_\infty$. By $K^\circ_\alpha$ we denote the subgroup
  in $K$ consisting of elements sending the set $\{ 1, \dots, \alpha\}$
  to itself. The quotient group $K^\circ_\alpha/K_\alpha$ is the finite symmetric
  group $S_\alpha$.

According
Arthur Lieberman \cite{Lie},
{\it any irreducible representation of the group $S_\infty$ is 
induced from an irreducible  representation of
some subgroup   $K^\circ_\alpha$ trivial on $K_\alpha$}.  
  
  \sm

  {\it The train of $S_\infty$.}
   Let $g\in S_{\infty}$. Denote by $(g)_{\alpha \beta}$
  the partial bijection $\{1, \dots, \beta\}\to\{1,\dots,\alpha\}$,
  which sends $i\le \beta$ to $j\le\alpha$ if and only if  $g$
  sends $i$ to $j$. It is easy to show that the map $g\mapsto(g)_{\alpha,\beta}$
  determines a bijection between
   double cosets $K_\alpha\backslash S_\infty/K_\beta$ and
   partial bijections $\{1, \dots, \beta\}\to\{1,\dots,\alpha\}$. 
   
   \sm

 According \cite{Olsh-sym} (see also exposition in \cite{Ner-book}, Sect. VIII.1-2),
  the following theorem
 holds:
 
 \begin{theorem}
 \label{th:lieb}
{\rm a)} The train 
$\Tra\bigl(S_\infty,\{K_\alpha\}\bigr)$
is the category whose objects are sets 
$\{1,\dots,\alpha\}$ and morphism are partial bijections.

\sm

{\rm b)} The functor $\rho\mapsto\wt \rho_{\alpha,\beta}$
is a one-to-one correspondence between equivalence classes
of unitary representations of $S_\infty$
and equivalence classes of $*$-representations of
$\Tra\bigl(S_\infty,\{K_\alpha\}\bigr)$.
 \end{theorem}
 
So unitary representations of $S_\infty$
(defined up to equivalence)
are in one-to-one correspondence with $*$-representations
of the category $\PB$ since this category is equivalent to $\Tra\bigl(S_\infty,\{K_\alpha\}\bigr)$.

 Classification of $*$-representations of the category 
 $\Tra\bigl(S_\infty,\{K_\alpha\}\bigr)$
 is a simple problem, see \cite{Ner-book}, Sect. VIII.1-2
 (also, see  below, Subsect. \ref{ss:rep-tra}).
 In fact, Theorem \ref{th:lieb} is equivalent to Lieberman's description \cite{Lie} of unitary representations of
 $S_\infty$.

% {\sc Remark.} This statement is equivalent to the Lieberman classification
% of unitary representations of $S(\infty)$ (any irreducible representation
% of $S(\infty)$ is induced (in the sence of Mackey \cite{Mach})
%  from an irreducible finite dimensional of the subgroup 
%  $S(n)\times S(\infty-n)$ (i.e., a group sending the set $\{1,2,\dots,n\}$
%  to itself, see, e.g, \cite{Olsh-sym},  \cite{Ner-book}, Sect. VIII.1-2.
%   It seems that representation theory of $S(\infty)$  contains few interesting
%   in comparision  with the representation theory of usual symmetric groups.
%   However the Lieberman theorem is an important point of a wider theory.
%    The existing 'representation theory of infinite symmetric group'
%(E.Thoma, A.M.Vershik, S.V.Kerov, G.I.Olshanski, A.Yu.Okunkov, A.M.Borodin} 
%is mainly representation theory of the following 'bysymmetric group'. We consider 
%group  consisting
%of pairs $(g_1,g_2)\in S(\infty)\times S(\infty)$ such that $g_1^{-1}g_2$
%is finitely supported,  the produt is $(g_1,g_2)(h_1,h_2)=(g_1h_1,g_2h_2)$.
%So we have a group $G$ containing $K=S(\infty)$, in this case we have
%a category of double cosets $K_\alpha \backslash G/K_\beta$, see
%\cite{Olsh-sym}. There are many other groups containing $S(\infty)$, which produce %numerous categories of  double cosets, see \cite{Ner-field}. 
%\hfill $\boxtimes$
  
 \section{Oligomorphic groups and their trains}
 
 \COUNTERS
 
 Here we construct categories of reduced double cosets for oligomorphic groups. This is a simple corollary
 of description of unitary representations of oligomorphic groups
 obtained by Tsankov \cite{Tsa}. These categories embed in  the semigroups of partial bijections
 examined by
Ben Yaacov, Ibarluc\'\i a, Tsankov \cite{BIT}. 

\sm 
 
  {\bf\punct Oligomorphic groups.}  A closed subgroup  $G$
  in an infinite  symmetric group $S(\Omega)$ is called {\it oligomorphic} if 
 it has a finite number of orbits on each multiple product 
 $\Omega^n:=\Omega\times\dots\times \Omega$ ($n$ times), see Cameron \cite{Cam}.
 See a collection of examples below in Subsect. \ref{ss:examples}.
 
\sm  
 
 The stabilizer of any finite subset in $\Omega$ is oligomorphic.
 
\sm  
 
If $H\subset G$ is an open subgroup, then $G$ is closed in the symmetric group
$S(G/H)$ and has a finite number of orbits on each set  $(G/H)^n$.

\sm

  Tsankov \cite{Tsa}
 obtained the following description of unitary representations of
 such groups:
 
 \begin{theorem}
 \label{th:tsankov}
 {\rm a)}  Each representation of
 $G$ is a direct sum of irreducible representations%
 \footnote{For this reason, the group $G$ has type $I$ 
 (see, e.g., \cite{Mack3}, Chaper 1).}.
 
 \sm
 
 {\rm b)} For any irreducible representation $\rho$ of $G$ there is an open subgroup
 $H^\circ$ and a normal subgroup $H$ in $H^\circ$ of finite index such that
 $\rho$ is induced\phantom{l}%
 \footnote{For induced representations of totally disconnected groups,
 see Mackey  \cite{Mack}, \cite{Mack3}. Also, see below Subsect. \ref{ss:tra-K}.}
 from a representation of $H^\circ$ trivial on $H$.
 \end{theorem}
 
 Irreducible representations induced from subgroups
 $H$ and $H'$ can be equivalent only for  trivial reasons:
 if $H$ and $H'$ are conjugate and their representations
 $\sigma$, $\sigma'$ are conjugate. 
 
\sm 
 
 {\bf \punct Perfect oligomorphic groups.}
 Fix an oligomorphic subgroup $K\subset S(\Omega)$.  For a  finite subset 
 $A\subset \Omega$ we
denote by $K(A)\subset S(\Omega)$ the subgroup consisting of all $g\in S(\Omega)$ fixing $A$ pointwise. We say that a finite subset
 $A\subset \Omega$ is {\it perfect}%
\footnote{In the terminology of model theory, 'a perfect set' is  'an algebraically closed set' 
(origin of this words is clear from Examples $3^\lozenge$, $4^\lozenge$
of Subsect. \ref{ss:examples} or Proposition \ref{pr:subspaces-perfect}),
see \cite{EH}, \cite{Tsa}. 
'Perfect oligomorphic groups' are 'groups of automorphisms of  first order
structures admitting a weak elimination of singularities', 
See also \cite{EH}, Lemma 1.3.}
 if all orbits of $K(A)$ on $\Omega\setminus A$ are infinite.   
  For a perfect $A$ denote by
  $K^\circ(A)\subset S(\Omega)$ the subgroup sending $A$
to itself. Denote
$$
\Gamma(A):=K^\circ(A)/K(A),
$$ 
it is a finite group acting on $A$.

\sm

%{\sc Example.} a) Let $S(\Omega)$ act on $\Omega$. Then all finite subsets are
%perfect. We have $K(A)=S(\Omega\setminus A)$, $K^\circ(A)=S(A)\times %S(\Omega\setminus A)$, $\Gamma(A)\simeq S(A)$. See more examples in Subsect. %\ref{ss:examples}.
%\hfill $\boxtimes$

\begin{lemma}
\label{l:perfect-set}
{\rm a)} Let $A\subset \Omega$ be perfect, $g\in K$. Then $gA$ is perfect.

\sm

{\rm b)} If $A$, $B$ are perfect, then $A\cap B$ is perfect.

\sm

{\rm c)} For any finite $R\subset \Omega$ there is a minimal
perfect set $\ov{\ov R}$ containing $R$.
\end{lemma} 

{\sc Proof.} Statements a), b) are obvious. 

Since $K$ is oligomorphic, the group $K(R)$ has 
a finite number of orbits on $\Omega$. In particular, the number of finite orbits is finite.
Let $\ov{\ov R}$ be the union of all finite orbits
of $K(R)$.  
 The group $K(\ov{\ov R})$ has finite index in $K(R)$ and 
elements of $\Omega\setminus \ov{\ov R}$, which have  finite orbits with respect to
$K(\ov{\ov R})$, have finite 
orbits with respect to $K(R)$. So such elements are absent.
\hfill $\square$

\sm

Notice that the set of  classes of $K$-equivalent perfect sets
 has a {\it natural partial order}: 
 $$B\preceq C$$
  if there is
$g\in K$ such that $gB\subset C$.

\sm 
 
We say that a closed subgroup $G$ in $S(\Omega)$  is {\it perfect oligomorphic}
  if it is oligomorphic and
  for each open subgroup $H\subset G$ there is a perfect finite subset
  $A$ such that $K(A)\subset H\subset K^\circ(A)$.
  
  \sm
  
 {\sc Remark.} This property depends on an embedding of $G$
in the infinite symmetric group.
  For instance, consider the set $\Omega^{[k]}$ of all $k$-element subsets
  in $\Omega$ ($k>1$) and the action of $S(\Omega)$ on $\Omega^{[k]}$. Clearly,
  the subgroup $S(\Omega)$ in $S\bigl(\Omega^{[k]}\bigr)$ is not perfect oligomorphic.
  \hfill $\boxtimes$

\sm

  Fix a perfect oligomorphic group $K\subset S(\Omega)$. Let $A$, $B$ be
  perfect subsets in $\Omega$. We say that a partial bijection $\theta:A\to B$
  is a {\it restriction} $g\bigr|_{A\to B}$ of $g\in K$ if  $\dom\theta=A\cap g^{-1} B$
  and the bijection $\dom \theta\to \im\theta$ is induced by $g$.
  
 \begin{lemma}
 \label{l:well}
 Let $K\subset S(\Omega)$ be perfect oligomorphic.
 Let $A$, $B$, $C\subset \Omega$ be perfect sets.
 Let  partial bijections  $\mu:A\to B$ , $\nu:B\to C$ have the form
 $$
 \mu=p\Bigr|_{A\to B},\quad \nu= q\Bigr|_{B\to C},
 \quad\text{where $q$, $p\in K$.}
 $$
 Then there is $r\in K$ such that $\nu\mu=r\Bigr|_{A\to C}$.
 \end{lemma}

We need the following P.~Neumann Lemma (see
\cite{Neu},  \cite{Cam}, (2.16)):

\begin{lemma}
\label{l:finite-shift}
 Let a group $G$ act on a countable set $\Xi$ without finite
orbits. Let $M$, $N$ be finite subsets in $\Xi$. Then there is 
$g\in G$ such that $gM\cap N=\varnothing$.
\end{lemma}
 
 {\sc Proof of Lemma \ref{l:well}.}
 We set 
$$\Xi:=\Omega\setminus B, \quad
G:=K(B),\quad
M:=p(A)\setminus B,\quad N:=q^{-1}C\setminus B
$$
and choose $g$ as in the Neumann Lemma. We set
$r=qgp$.
\hfill $\square$

\sm
 
{\bf\punct Trains of perfect oligomorphic groups.}
  Consider the following
  category $\Tra(K)$ (the {\it train of a perfect oligomorphic group $K\subset S(A)$}). Its objects are perfect finite subsets in $\Omega$.
 A {\it morphism} $\theta:A\to B$ is a partial bijection
  that can be obtained by a restriction of $g\in K$ to $A$.
  A product of morphisms is the product of partial bijections
  (by Lemma \ref{l:well} such product is a morphism). The involution
  in $\Tra(K)$ is the involution in the category of partial bijections.
  
  \sm
  
{\sc Remark.} So $\Tra(K)$  is  realized by partial bijections
and the set of morphisms is closed with respect to the involution.
Therefore, it is an 'inverse category' in the sense of \cite{Kas}, i.e., a categorical
counterpart  of inverse semigroups. %, see, e.g., \cite{CP}.
\hfill $\boxtimes$
  
\begin{theorem}
\label{th:main}
Let $K\subset S(\Omega)$ be a perfect oligomorphic group, let $A$ range 
in the space of perfect subsets in $\Omega$.

\sm

{\rm a)} The category $\Tra(K)$ is the train $\Tra(K,\{K(A)\})$ of
the group $K$.

\sm

{\rm b)} The map $\rho\mapsto \wt \rho$ is a bijection
between equivalence classes of unitary representations  of $K$ and
equivalence classes of $*$-representations of the
category $\Tra(K)$.
\end{theorem} 

{\sc Remark.} We also can reduce $\Tra(K)$ choosing one representative from 
each $K$-equivalence class  of perfect sets. In this way, we pass
to an equivalent category, which is 'smaller' and more transparent.
\hfill $\boxtimes$

\sm

The proof of this statement given below in 
Subsect. \ref{ss:tra-K}--\ref{ss:rep-tra} is a simple corollary
of Tsankov's Theorem \ref{th:tsankov}, formal manipulations with 
'ordered categories' in the sense of \cite{Ner-book}, and Mackey's
technique of induced representations \cite{Mack}, \cite{Mack2}, 
\cite{Mack3}, \cite{BR}.

\sm

{\bf\punct Description of representations of $\Tra(K)$ corresponding
to representations of $K$.%
\label{ss:tra-K}} Here we verify
the statement a) of Theorem \ref{th:main}.

According  Theorem \ref{th:tsankov}, any irreducible representation
of a perfect oligomorphic group $K$ is induced by
 an irreducible representation of some subgroup
$K^\circ(A)$ trivial on $K(A)$. 

\sm

 {\it Let us reformulate a definition
of induced representations} \cite{Mack} {\it in a convenient for us form.}
Let $A\subset \Omega$ be a perfect set. Denote by $\cE(A)$ the set of all maps
$\iota: A\to \Omega$ obtained by restrictions of elements $g\in K$. Clearly,
$\cE(A)$ is a $K$-homogeneous space $K/K(A)$, the group $K$
acts on $\cE(A)$ by  left multiplications $\iota\mapsto g^{-1}\iota$.
We also have an action of the finite group $\Gamma(A)=K^\circ(A)/K(A)$
 on $\cE(A)$ by right
multiplications $\iota\mapsto \iota\gamma$. These two actions commute.

Consider a unitary irreducible representation $\sigma$ of  $\Gamma(A)$
in a Euclidean space $V$. Consider the space
$\ell^2(\cE(A), V)$ of $V$-valued functions on $\cE(A)$ equipped with
the $\ell^2$-inner product
$$
\la F_1, F_2\ra_{\ell^2}=
\sum_{\iota\in \cE\,(A)} \la F_1(\iota), F_2(\iota)\ra_V.
$$
The group $K$ acts in  $\ell^2(\cE(A), V)$ by shifts
\begin{equation}
\pi(g)
F(\iota)=F(g^{-1}\iota).
\label{eq:shift}
\end{equation}
 We consider the subspace 
$H=H(A,\sigma)\subset\ell^2(\cE(A), V)$ consisting of functions
satisfying
\begin{equation}
F(\iota\gamma)=\sigma(\gamma^{-1}) F(\iota),
\label{eq:iota-gamma}
\end{equation}
where $\gamma\in\Gamma(K)$.
Restricting operators \eqref{eq:shift} to $H(A,\sigma)$, we get a
unitary representation 
\begin{equation}
\tau=\Ind_{K^\circ(A)}^K(\sigma)
\label{eq:Ind}
\end{equation}
 of $K$ in $H(A,\sigma)$
(the representation of the group $K$ {\it induced} by the representation 
$\sigma$ of the subgroup $K^\circ(A)$).

\sm

{\sc Remark.} So, the representation $\tau$
is a subrepresentation 
in the quasiregular representation of $K$ in
$
\ell^2\bigl(\cE(A)\bigr)=\ell^2\bigl( K/K(A) \bigr)
$. \hfill $\boxtimes$

\sm

Let us examine the operators corresponding to double cosets.
For a perfect set  $B$, denote by $\cE(A,B)$ the set of all embeddings
$\iota:A\to B$ such that $\iota=g\bigr|_{A}$ for some $g\in K$, so $\cE(A,B)\subset \cE(A)$.
Denote by $H(A,\sigma)_B\subset H$ the subspace of $K(B)$-invariant functions, by $P(B)$ the projection to 
$H(A,\sigma)_B\subset H(A,\sigma)$.

\begin{lemma}
The subspace $H(A,\sigma)_B$ consists of functions $F\in H(A,\sigma)$
supported by $\cE(A,B)$. The projection $P(B)$
is the operator of
restriction of  functions from $\cE(A)$ to the subset $\cE(A,B)$.
%
% multiplication by the characteristic function of the set
%$\cE(A,B)$
\end{lemma}
  
{\sc Proof.} A $K(B)$-invariant function $F\in \ell^2$ must be constant on
$K(B)$-orbits. If $\iota (A)\not\subset B$, then $K(B)$-orbit of $\iota$
is infinite. Since $F\in \ell^2$, it is zero on such orbits.

If $\iota(A)\subset B$, then $\iota$ is fixed by the subgroup $K(B)$.
Therefore functions $F$ supported by $\cE(A,B)$ are fixed. 
\hfill $\square$

\sm

Fix perfect sets $B$, $C$. For any $g\in K$ we assign 
the partial bijection $g\Bigr|_{C\to B}$ as above. Clearly,
the map $g\mapsto g\Bigr|_{C\to B}$ is constant on double cosets
$K(B)\backslash K/K(C)$. So, we have a surjective (generally, non-injective)
map 
$$K(B)\backslash K/K(C)\to \Mor(C,B).$$

Now we define
a representation $(T,\tau)$ of the category $\Tra(K)$.
A space $T(B)$ is the space of $V$-valued functions on $\cE(A,B)$
satisfying \eqref{eq:iota-gamma}.
An operator
$$
\tau_{B,C}(g)=P(B)\, \tau(g)\Bigr|_{T(C)}: T(C)\to T(B) 
$$
is given by the formula:
$$
\tau_{B,C}(g)F(\iota)=\begin{cases}
F(g^{-1}\iota ),\quad&\text{if $i(A)\subset g^{-1}C\cap  B$};\\
0, \quad&\text{otherwise}.
\end{cases}
$$

Clearly, we get a representation of the category $\Tra(K)$
and this establishes
 the statement a) of Theorem \ref{th:main}.
We separately formulate the following obvious statement.

\begin{lemma}
\label{l:EndAut0}
For the representation $(T,\tau)$ of $\Tra(K)$, a space $T(B)$
is nonzero if and only if $B\succeq A$. The representation
of $\End(A)$ in $T(A)$ is zero on $\End(A)\setminus \Aut(A)$
and coincides with $\sigma$ on $\Aut(A)$.
\end{lemma}

\sm

{\bf \punct Representations of the category $\Tra(K)$.
\label{ss:rep-tra}}
It remains to verify that each irreducible $*$-representation
of the category $\Tra(K)$ arises from an irreducible representation of $K$.

%Notice that the set of equivalence claseses of objects of $\Tra(K)$
% has a natural partial order: $B\prec C$ if there is an embedding
% $\mu:B\to C$ (it satisfies $\mu^* \mu=1_B$). 

\begin{lemma}
\label{l:EndAut}
{\rm a)} Let $(Q,\kappa)$ be a $*$-representation of $\Tra(K)$. Let $D$
 be a $\preceq$-minimal perfect set such that $Q(D)\ne 0$.
 Then $\kappa_{\fan D,D}(\phi)=0$ for any morphism $\phi\in \End(D)\setminus \Aut(D)$.
 
 \sm 

{\rm b)} Additionally, let  $(Q,\kappa)$ be  irreducible. Then a
 $\preceq$-minimal perfect set  such that 
$Q(D)\ne 0$ is unique {\rm(}upto $K$-equivalence{\rm)}. 
\end{lemma}

We call this $\preceq$-minimal perfect set $D$ by the  {\it conductor}
of an irreducible representation.

\sm

{\sc Proof.} Let $\phi\in \End(D)\setminus \Aut(D)$, $\kappa_{D,D}(\phi)\ne 0$.
Then the partial bijection $\phi^*\phi$ is an identical map
  $\dom \phi\to\dom\phi$, and 
$$\kappa_{\fan D,D}(\phi^*\phi)=\kappa_{\fan D,D}(\phi)^*\kappa_{\fan D,D}(\phi)\ne 0
.$$
Indeed, if $A:L_1\to L_2$ is a non-zero operator between Hilbert spaces,
then $A^*A$ also is non-zero.

 By Lemma \ref{l:perfect-set}.a)-b),
$\dom \phi$ is a perfect set, denote it by $C$. Let $\mu$ be the identical embedding
$C\to D$, so $\mu\mu^*=\phi^*\phi$. Therefore, $\kappa_{\fan  D,C}(\mu)\ne 0$.
Hence,  
$$\kappa_{\fan  C,C}(\mu^*\mu)=
\kappa_{\fan  D,C}(\mu)^*\kappa_{\fan  D,C}(\mu)\ne 0.$$ 
But $\mu^*\mu$ is the unit in $\End(C)$,
and therefore $Q(C)\ne 0$. We come to a contradiction.

\sm

b) Let $D$, $D'$ be two non-equivalent $\preceq$-minimal perfect sets. 
Let they be non-equivalent. Let $\psi\in \Mor(P,P')$.
If $\dom \psi \ne P$, then $\psi^* \psi\in \End(P)\setminus\Aut (P)$.
So, $\kappa(\psi^* \psi)\ne 0$. Hence $\kappa(\psi)=0$.

If $\im \psi\ne P'$, then $\psi\psi^*\in \End(P')\setminus\Aut (P')$.
Again, $\kappa(\psi)=0$. 

So, $\kappa(\psi)=0$ for all $\psi\in \Mor(P,P')$.
This contradicts to irreducibility.
\hfill $\square$ 

\sm

Now Theorem \ref{th:main} follows from straightforward arguments,
which are presented below.

Let $(Q,\kappa)$ be a representation of a category $\cK$. Let $C\in \Ob(\cK)$.
 Let
$\xi_\alpha$ be a family of vectors in $Q(C)$.
For any $D\in\Ob(\cK)$ we consider the closed linear
span $N(D)$ of all vectors $\kappa_{\fan D,C}(\mu)\xi_\alpha$, where $\mu$ ranges in $\Mor(C,D)$. Clearly, the family of subspaces $N(D)$ is a subrepresentation in $(Q,\kappa)$
(we call it the {\it cyclic span} of $\xi_\alpha$). We say that a family $\xi_\alpha$ is 
{\it cyclic} if $N(D)=Q(D)$ for all $D$.

\begin{lemma}
\label{l:GNS}
If for a  $*$-representation $(Q,\kappa)$ of  $\Tra(K)$ a  subspace $Q(C)$ is cyclic, then $(Q,\kappa)$
is uniquely determined
by the representation of $\End(C)$ in $Q(C)$.
\end{lemma}

{\sc Proof.} This is a version of the GNS-construction.
Let $v_1$, $v_2\in Q(C)$, $r_1$, $r_2\in \Mor(C,D)$.
Then
\begin{multline*}
\bigl\la \kappa_{\fan D,C} (r_1)v_1,\, \kappa_{\fan D,C} (r_2)v_2\bigr\ra_{Q(D)}=
\bigl\la \kappa_{\fan D,C} (r_2)^* \kappa_{\fan D,C} (r_1)v_1,\, v_2\bigr\ra_{Q(C)}
=\\=\bigl\la \kappa_{\fan C,C} (r_2^* r_1)v_1,\, v_2\bigr\ra_{Q(C)},
\end{multline*}
and we see that such inner products are uniquely determined by
the representation of $\End(C)$. Let $(Q',\kappa')$ be another representation
of the category $\Tra(K)$, which coincides with $(Q,\kappa)$ on $\End(C)$.
Then for each $D$ we can identify $Q(D)$ with $Q'(D)$ assuming 
that any vector $\kappa'_{\fan D,C}(r) w$ corresponds to 
$\kappa_{\fan D,C}(r) w$.

This family of operators $Q(D)\to Q'(D)$ is an intertwiner. 
Indeed, let $q\in \Mor(D,E)$, $w\in Q(C)$. Let us evaluate  
inner products of $\kappa_{\fan E,D}(q)\kappa_{\fan D,C}(r)w$ with vectors 
$\in Q(E)$:
$$
\bigl\la \kappa_{\fan E,D}(q)\kappa_{\fan D,C}(r)w,\, \kappa_{\fan E,C}(p) v\bigr\ra_{\fan Q(E)}
=\bigl\la \kappa_{\fan C,C}(p^*qr)w,\, v \bigr\ra_{Q(C)}.
$$
So, these inner products are uniquely determined by
the representation of the semigroup $\End(C)$.
Hence, 
 a position  of $\kappa_{\fan E,D}(q)\kappa_{\fan D,C}(r)w$ also is determined by
 the representation
 $\kappa_{\fan C,C}$.
 \hfill
 $\square$
 
\sm

\begin{lemma}
Any $*$-representation $(R,\rho)$ of $\Tra(K)$ is a direct sum of irreducible representations. 
\end{lemma} 

{\sc Proof.} 
 Consider a sequence $A_1$, $A_2$, $\dots$ including all perfect sets
 in $\Omega$ defined up to the $K$-equivalence. Since $K$ is oligomorphic,
 it has only finite number of orbits on the space of finite
 subsets of a fixed cardinality. So, we can assume
that numbers of elements of $A_j$ are (non-strictly) increasing.
We take the cyclic span $(Q^{(1)},\kappa^{(1)})$ of $R(A_1)$ and the orthogonal complement 
$(R^{(1)},\rho^{(1)})$ of the cyclic span. By  construction, 
we have $R^{(1)}(A_1)=0$. Next, we take the cyclic span $(Q^{(2)},\kappa^{(2)})$
of $R^{(1)}(A_1)$ and  orthogonal complement $(R^{(2)},\rho^{(2)})$ of $(Q^{(2)},\kappa^{(2)})$ in $(R^{(1)},\rho^{(1)})$.
 Now we have
$R^{(2)}(A_1)=0$, $R^{(2)}(A_2)=0$. Etc. 

In this way, we get a decomposition $(R,\rho)=\oplus_j (Q^{(j)},\kappa^{(j)})$.
In each summand $(Q^{(j)},\kappa^{(j)})$, a perfect set $A_j$
is a minimal set for which $Q^{(j)}(\cdot)\ne 0$. By Lemma \ref{l:EndAut}, the  action
of $\End(A_j)$  in $Q^{(j)}(A_j)$ is reduced to the action of $\Aut(A_j)$.
We decompose each representation of a finite group $\Aut(A_j)$
into a direct sum, take the cyclic span  of each summand with respect to $\Tra(K)$, and come to a decomposition of each $(Q^{(j)},\kappa^{(j)})$ into a direct sum of irreducible representations. 
\hfill
$\square$

\sm 

According  Lemmas \ref{l:EndAut} and \ref{l:GNS},
an irreducible $*$-representations
$(R,\rho)$
of the category $\Tra(K)$ is uniquely determined by 
its conductor $A$   and 
an irreducible representation $\sigma$ of a finite group $\Aut(A)$.

According Lemma \ref{l:EndAut0}, these representations coincide with 
representations originated from representations of the group $K$.
This completes the proof of Theorem \ref{th:main}.

\sm

%{\sc Remark.} In the proof above, we actually used that our category:
%
%\sm 
%
%---  is realized as 
% a subcategory of $\PB$ and is closed with respect to the involutions (taking of a pseudoinverse);
%
%\sm 
%
%--- for any subset $\dom \mu$, where $\mu$ is a morphism, there is a bijective  partial bijection $\in \Mor(A,\dom\mu)$
%for some object $A$. 
%
%\sm 
%
%Under these conditions an irreducible representations of the category is uniquely
%determined by its conductor $D$  and an irreducible representation of the group $\Aut(D)$.
% \hfill $\boxtimes$

\sm   

{\bf\punct Weak closures.%
\label{ss:weak}} Let $\rho$ be a unitary representation of 
a group $G$. Consider the image $\rho(G)$ in the unitary group
and its   closure $\ov{\rho(G)}$ in the space all operators with respect to the weak operator
 topology.
It is well-known (and can be easily verified, see, for example, \cite{Olsh-semigroup},
  \cite{Ner-book}, Sect. I.4)
that $\ov{\rho(G)}$
 is a compact metrizable semigroup with a separately continuous multiplication.
 See a discussion of such closures for 'infinite-dimensional' ('large') groups in 
 \cite{Olsh-semigroup}, \cite{Ner-book}, for oligomorphic groups in 
 \cite{BIT}.

\begin{proposition}
Let $K\subset S(\Omega)$ be a perfect oligomorphic group.
Let $B\subset \Omega$ be a perfect set.
Then there is a sequence $g_j\in K(B)$ such that
  for any unitary representation
$\rho$
of $K$
 the orthogonal projection $P(B)$ to the subspace of $K(B)$-fixed vectors
is a weak operator limit of the sequence $\rho(g_j)$.
\end{proposition}

{\sc Remark.} As far as the author knows, a similar statement holds in all cases,
 when
multiplication of double cosets exists.
\hfill $\boxtimes$

\sm

{\sc Proof.} It is sufficient to verify the statement for irreducible representations.
 Consider an exhaustive family of subsets 
$B\subset R_1\subset R_2\subset\dots\subset \Omega\setminus B$. By Lemma \ref{l:finite-shift},
for each $j$ there exists $g_j\in K(B)$ such that 
$$g_j(R_j\setminus B)\cap (R_j\setminus B)=\varnothing.$$

Consider the quasiregular representation $\pi$ of $K$ in $\ell^2(\cE(A))$, see \eqref{eq:shift}. 
For $\iota\in \cE(A)$ denote by $e_\iota$ the $\delta$-function supported 
by $\iota$. If $\iota\in\cE(A,B)$, then $\iota$ is fixed by all $g_j$.
Fix $\iota(A)\not\subset B$. Let $\nu\in \cE(A)$.
There is $N$ such that $(\iota(A)\setminus B)$, $(\nu(A)\setminus B)\subset R_N$.
For large $j\ge N$ we have 
$$\bigl(g_j(\iota(A)\setminus B)\bigr) \cap \bigl(\nu(A)\setminus B\bigr)=\varnothing.$$
Therefore,  $g_j\iota\ne \nu$. Hence,
$$
\la\pi(g_j) e_\iota, e_\nu\ra_{\ell^2}=
\la e_{g_j\iota},e_\nu\ra_{\ell^2}=\delta_{g_j\iota, \nu}=0.
$$
Therefore, the weak limit of the sequence $\pi(g_j)$ is
the projection to the space  of $K(B)$-fixed vectors.

Any induced representation \eqref{eq:Ind} is a subrepresentation of
$\pi$, and we get the desired statement.
\hfill 
$\square$ 

\sm

\begin{corollary}
For any unitary representation $\rho$ of a perfect oligomorphic group 
$K\subset \Omega$
for each perfect sets $A$, $B$
operators 
$$
\wt \rho_{A,B}(g)=P(A) \,\rho(g) P(B)
$$
are contained in the weak closure of $\rho(K)$.
\end{corollary}

\sm

So, the category $\Tra(K)$ has a canonical embedding in the universal
universal  semigroup compactification of $K$,
which was examined by Ben Yaacov, Ibarluc\'\i a, Tsankov \cite{BIT}.

\sm 

{\bf\punct Examples.%
\label{ss:examples}} Here we discuss an initial collection of examples
of
oligomorphic groups (for $1^\lozenge$-$4^\lozenge$,
 see \cite{Tsa}, Sect. 5), and the corresponding categories.

\sm

$0^\lozenge.$
Consider the infinite symmetric group $S(\N)$, i.e., the group of all permutations of
$\N$. All finite subsets in $\N$ are perfect,
morphisms are partial bijections. This category is equivalent
to the category $\PB$, 
see Subsect. \ref{ss:S-infty}.

\sm

$1^\lozenge.$ Consider the group $G=\Aut(\Q,<)$ of all order preserving bijections of the set
$\Q$ of rational numbers%
\footnote{$\Q$ is a 'generic' contable linearly ordered set.
For instance, let us construct a random ordered set in the following way.
 We take a point $x_1$. 
With probability $1/2$ we assume $x_2<x_1$ or $x_2>x_1$. With probability $1/3$ we put $x_3$ to one of possible 3 positions, etc. With probability 1 we come to  an ordered set equivalent to $\Q$.}. Perfect sets are arbitrary finite subsets  $A\subset\Q$.
We come to the category  $\Tra\bigl(\Aut(\Q,<)\bigr)$. It objects are
finite  subsets of $\Q$. Morphisms $A\to B$ are order preserving partial bijections.
Groups of automorphisms are trivial.

This category is equivalent to the category of finite 
ordered sets, morphisms are order preserving partial bijections.
Another equivalent category is the category, whose objects
are sets $\{1, \dots,\alpha\}$, where $\alpha=0$, 1, 2, \dots;
morphisms are order preserving partial bijections.   

So, irreducible unitary representations
are enumerated by the  number $\# A$ of elements of a conductor $A$.
In other words, any irreducible representations of this group is 
realized in the spaces $\ell^2$ on the set of all $\alpha$-element subsets
in $\Q$ for some $\alpha$.

\sm

$2^\lozenge.$ Let $\cR$ be the Rado graph,  see, e.g., \cite{Cam}.
Recall that it is a graph with a countable number of vertices satisfying
the following  condition (the 'Urysohn extension property'): let $A\subset B$ be finite graphs,
then any embedding $A\to\cR$ extends to an embedding $B\to\cR$. 

Consider 
 a countable set of vertices, for each pair of vertices we draw an
edge connecting them with probability $1/2$. Then we get the Rado graph with probability 1, see \cite{ER}, \cite{Cam-Rad}.

 The group $\Aut(\cR)$ of automorphisms of $\cR$ acting on the set
of vertices is perfect oligomorphic, arbitrary finite subset of the set of vertices is perfect. Objects of the train $\Tra(\Aut(\cR))$
are finite  subgraphs $X\subset\cR$, a morphism $X\to Y$ is an isomorphism
of an induced subgraph%
\footnote{To avoid a misleading we use the term '{\it induced subgraph}': we consider some set of vertices of
	$X$ and {\it all} edges of $X$ connecting these vertices.} $A\subset X$ to an induced subgraph $B\subset Y$.
	
	This category is equivalent to the category, whose
	objects are finite graphs and morphisms are isomorphisms
	of induced subgraphs.

Irreducible representations are enumerated by pairs
\begin{multline*}\Bigl(\text{a finite graph $\Delta$ defined up to   isomorphism,}
	\\
	\text{an irreducible representation of the group $\Aut(\Delta)$  defined
		up to
		equivalence}
	\Bigr).
\end{multline*}

Recall that any finite group can be realized as a group of automorphisms
 of some finite graph.

\sm

$3^\lozenge.$    Consider a linear space $\ell(\F)$ over
a finite field $\F$ consisting of sequences $v=(v_1,v_2,\dots)$,
where $v_j\in\F$ and $v_k=0$ for sufficiently large $k$. 
The group $\GL(\infty,\F)$ of all linear transformations
of $\ell(\F)$ can be regarded as the group of all invertible infinite matrices
 over $\F$,
which have finite number of elements in each column. 
The group $\GL(\infty,\F)$ is perfect oligomorphic,
perfect subsets are finite-dimensional linear subspaces
in $\ell(\F)$.

Objects of the train are finite dimensional spaces $X$ in $\ell(\F)$.
A morphism $X\to Y$ is a {\it partial linear bijection}, i.e.,
 an invertible operator from a subspace  $A\subset X$ to a subspace
  $B\subset Y$. Groups of automorphisms are group of linear transformations
  $\GL(X)$.
  
This category is equivalent to the category of all
finite-dimensional linear spaces over $\F$ (with same morphisms).
We also can consider another  equivalent category, whose objects
are coordinate spaces $\F^0$, $\F^1$, $\F^2$, \dots with same morphisms.   
  
   So irreducible  representation are enumerated by
  a number $n=0$, $1$, $2$, $\dots$ and a complex irreducible
  representation of the group $\GL(n,\F)$; the picture
  is completely parallel to representations of $S_\infty$. 
  
  \sm
  
  $4^\lozenge.$ (\cite{Olsh-semigroup}, \cite{Tsa})
Consider the group $\GL^\diamond(\infty,\F)$ consisting of infinite invertible matrices over  a finite  field $\F$, which have finite number of elements in each row and each column (inverse matrices must have the same form). Let $\ell(\F)$ be as in the previous example, let 
$\ell'(\F)$ be a copy of $V$.
Let 
$\GL^\diamond(\infty,\F)$ act in $\ell(\F)\sqcup \ell'(\F)$ by the formula
$v\mapsto gv$, $v'\mapsto g^{t-1}v'$, where $v\in \ell(\F)$,
$v'\in\ell'(\F)$, and $^t$ denotes the transposition.
Notice that this action preserves the natural pairing $\ell(\F)\times \ell'(\F)\to \F$
defined by
$
\{v,v'\}=\sum v_jv'_j
$.
The group $\GL^\diamond(\infty,\F)$ is a perfect oligomorphic
subgroup  in the symmetric group $S(\ell(\F) \sqcup \ell'(\F))$, any perfect subset in $\ell(\F)\sqcup \ell'(\F)$
is a union  $X\sqcup X'$ of a subspace  
$X\subset \ell(\F)$ and a subspace
 $X'\subset \ell'(\F)$. Morphisms $(X,X')\to (Y,Y')$
are pairs of partial linear bijections $\mu:X\to Y'$, $\nu:X'\to Y'$
such that for any $v\in\dom\mu$, $v'\in\dom\nu$ we have
$\{\mu(v),\nu(v')\}=\{v,v'\}$. 

We can consider the following equivalent category.
Objects are  pairs $(X,X')$ of finite-dimensional linear spaces
over $\F$ equipped with  a bilinear form $\Phi :X\times X'\to\F$.
Morphisms $(X,X';\Phi)\to (Y,Y';\Psi)$ are pairs of linear partial bijections $\lambda:X\to Y$ and $\mu: X'\to Y'$
such that $\Psi(\lambda(x),\mu(x'))=\Phi(x,x')$. 

 We consider the subgroup
$\GL(X\Join X',\Phi) $ in $\GL(X)\times \GL(X')$ preserving $\Phi$.
Irreducible representations of $\GL^\diamond(\infty,\F)$ (and of its train)  are enumerated by the following data (up to the natural equivalence):
 a pair $(X,X')$ equipped  with $\Phi$ and an irreducible representation of the group $\GL(X\Join X',\Phi) $.

Notice that pairs $(X,X';\Phi)$ are classified
(up to a natural equivalence) by $\xi=\dim X$, $\xi'=\dim X'$, and
the rank $\rho$ of $\Phi$. The group $\GL(X\Join X',\Phi)$
is a semidirect product
$$
 \Bigl(\GL(\xi-\rho,\F)\times \GL(\xi'-\rho,\F)
\times \GL(\rho,\F)\Bigr)\ltimes
\Bigl( (\F^\rho\otimes \F^{\xi-\rho})\times 
(\F^{\xi'-\rho}\otimes \F^{\rho})\Bigr). 
$$
For a classification of irreducible representations of
$\GL(X\Join X',\Phi)$ we can apply Mackey's arguments
\cite{Mack-ext} but they lead to a branching recursive problems
and do not give a final answer.

\sm

$5^\lozenge.$ Let $\F$ be a finite field of characteristic $\ne 2$.
Consider a non-degenerate  quadratic form $\Phi$ on $\ell(\F)$, all such forms are
equivalent%
\footnote{Let  $\epsilon\in \F$ be a quadratic non-residue. We reduce a quadratic form to a sum of squares
$\sum a_j w_j^2 $, where $a_j$=1 or $\epsilon$, in the usual
way. Binary forms $z_1^2+z_2^2$ and $\epsilon z_1^2+\epsilon z_2^2$
are equivalent. Therefore, each  form is equivalent to the form
$\sum w_{2j}^2+\sum \epsilon w_{2j+1}^2$.},
 so we can assume $\Phi(v)=\sum v_j^2$.
Consider   the group $\mathrm{O}(\infty)$ of linear transformations
 preserving this form.
It is perfect oligomorphic on $\ell(\F)$, perfect subsets are
finite-dimensional linear subspaces in $\ell(\F)$ (the proof for cubic forms
from the next section is valid for quadratic forms).
 A morphism $X\to Y$ is
a bijective linear map $\nu$ from a subspace $A\subset X$ to a subspace
$B\subset Y$ such that $\Phi(\nu(x_1),\nu(x_2))=\Phi(x_1,x_2)$.

We can also consider the following  equivalent category. Its objects are 
finite dimensional linear spaces $V$ over $\F$ equipped with quadratic forms $\Phi_V$, morphisms are partial linear bijections preserving 
forms. Irreducible unitary representation
are enumerated by quadratic spaces $(V,\Phi_V)$ defined up to the equivalence
(i.e., by $\dim V$, rank of $\Phi_V$, and the discriminant of $\Phi_V$)
and by irreducible representations of groups
$\rm{O}(V,\Phi_V)$ of isometries of $(V,\Phi_V)$.

%Notice that the latter group is a semidirect product
%of
%$$\GL(\rk \lambda, \F)\times \GL(\# X-\rk \lambda,\F)\times \GL(\# X'-\rk \lambda,\F)$$
%and an Abelian normal subgroup. 

 % \sm 

%In the next section we produce more examples in this spirit.

%{\sc Remark.} In all our examples except  $G=S_\infty$
%double cosets  really are reduced.
%\hfill $\boxtimes$  

\sm

{\bf \punct Some remarks and questions.\label{ss:question}}
1) In  Examples $1^\lozenge$, $2^\lozenge$, $4^\lozenge$, $5^\lozenge$,
 we  really have reduced double cosets. 

For the group $K=\Aut(\Q,<)$, double cosets $K(A)\backslash K/K(B)$
are enumerated by the following data. We consider an ordered set $Q$
and two order preserving injective maps $\lambda:A\to Q$, $\mu:B\to Q$
such that $\lambda(A)\cup \mu(B)=Q$. Taking the product $\lambda^*\mu$
in the sense of partial bijections we get 
the corresponding element of the reduced set 
$ [\![K(A)\backslash K/K(B)]\!]$.

For the group $\Aut(\cR)$, picture is similar.
Let $A$, $B$ be finite induced subgraphs of $\cR$.
Double cosets $K(A)\backslash K/K(B)$ are enumerated by the following data:

\sm

--- a finite graph $Q$;

\sm

--- embeddings $\lambda:A\to Q$, $\mu:B\to Q$ establishing 
isomorphisms of $A$ and $B$ with induced subgraphs in $Q$
such that $\lambda(A)$ and $\mu(B)$ covers the set of vertices
of $Q$.

\sm

In both cases, I do not see a possibility to define a natural operation
$$
K(A)\backslash K/K(B)\,\times\, K(B)\backslash K/K(C)\,\to\,
K(A)\backslash K/K(C),
$$ 
the multiplication arises after reduction.

For the group $\GL^\diamond(\infty,\F)$
it is possible to assign a product for any two double cosets
(it is a straightforward copy of a similar operation for real
classical groups). This product  was examined and clarified  by Olshanski  \cite{Olsh-semigroup}.
But in this case unitary representations do not separate double cosets,
and we have a strong reduction.

\sm

2) As we mentioned in  Introduction, theory of oligomorphic groups arose not in representation 
theory (see \cite{Cam}, \cite{Tsa}, \cite{BIT}).
From the point of view of refined representation theory
(as semisimple Lie groups, symmetric groups, Chevalley groups, 
$p$-adic groups,
and  in\-fi\-nite-di\-men\-si\-onal groups) their unitary representations can
seem exotic 
from one hand and
relatively  poor from the other hand.

We have seen a classifications of irreducible unitary representations 
of the group $\Aut(\cR)$, which includes classification of all irreducible representations of all finite groups. 
In this sense, the situation looks  very complicated and very rich. Nevertheless, visible to the author problems of representation 
theory for the group $\Aut(\cR)$ are reduced using the Mackey machinery
of induced representations \cite{Mack2}, \cite{Mack3}, \cite{BR}  to problems for  finite groups 
(which in this generality are unsolvable). Certainly, it is nice and impressive that
the group $\Aut(\cR)$ admits a classification of unitary representations. But it is not clear how to continue this,
if we want a continuation in the representation-theoretical spirit
(but   see \cite{Tsa1} on actions on measure spaces and de Finetti
type theorems). 

\sm 

3) At  first glance, we have a similar picture for complete infinite symmetric group $S_\infty$ and the groups $\GL(\infty,\F)$.
Their irreducible representations are enumerated by 
all irreducible representations of all  finite symmetric groups $S_n$
and all finite  groups $\GL(n,\F)$ respectively.
%$\GL^\diamond(\infty,\F)$.
Natural representation-theoretical problems for $S_\infty$ 
and  $\GL(\infty,\F)$  are reduced
to heavy representation-theoretical problems  for finite symmetric groups
$S_n$ and finite groups $\GL(n,\F)$. It seems that  we do not meet new entities
passing from $S_n$ to $S_\infty$ and from $\GL(n,\F)$ to $\GL(\infty,\F)$. 

% According Arthur Lieberman's
%classification \cite{Lie},  its irreducible unitary
%representations are enumerated by a pair (a nonnegative integer $n$,
%an irreducible representation of the symmetric group $S_n$).  Again, problems of 
%representation theory of $S_\infty$ are reduced to similar problems for finite  groups
%$S_n$ (which are classical and heavy). 

However, Lieberman classification  theorem (see \cite{Lie}, \cite{Olsh-sym}, \cite{Ner-book}) for representations of
$S_\infty$ is an important element of a wider 'representation
theory of infinite symmetric groups' (see, for example, \cite{Olsh-sym}, \cite{KOV}, \cite{Ner-field}, \cite{BO}). %\cite{Str}).
 This theory
 was initiated by the work of Elmar Thoma \cite{Tho}
(1964) on 'characters of infinite symmetric group'. Olshanski \cite{Olsh-sym} explained that  the Thoma characters correspond 
to unitary representations of the following '{\it bisymmetric group}': we consider
the subgroup in $S_\infty\times S_\infty$ consisting
of pairs $(g_1,g_2)$ such that $g_1g_2^{-1}$ are finitely supported
permutations%
\footnote{Consider an irreducible unitary representation $\rho$
	of a compact group $K$ in a Euclidean space $V$. Consider the space
	$\Op(V)$
	of all linear operators in $V$ equipped with the inner product $\la A,B\ra:=\tr AB^*$.
	The group $K\times K$ acts in $\Op(V)$ by operators
	$R(k_1,k_2)\,A= \rho(k_1)A\rho(k_2)^{-1}$. The identity operator $1_V$ is
	a unique vector fixed by the diagonal subgroup $\diag(K)\subset K\times K$.
	The matrix element of $1_V$ is
	$\la R(k_1 k_2) 1_V,1_V\ra=\tr \rho(k_1k_2^{-1})$. It is a value of the character of $\rho$
	at $(k_1k_2^{-1})$. So (in the standard terminology) 'characters of a group $K$ are 
	functions on $K\times K$ spherical with respect to $\diag(K)$.
	 %	Recall that characters of a finite or compact group $K$ are precisely
%functions on $K\times K$ that are spherical with respect to the diagonal subgroup. 
It turns out that a pass
to infinite limit in theories of characters (and, more generally, in
theories of  spherical functions)
is more interesting than attempts of straightforward generalizations of representations ('bisymmetric group' arises in such limit). Also,
it turns out that theory of Thoma characters can be included to a wider picture that (at least now) 
has no counterpart for finite symmetric groups, see \cite{Ner-field}.}. This group contains the diagonal consisting
of pairs $(g,g)$, where $g\in S_\infty$.  There is a big zoo
of other interesting groups containing $S_\infty$, see \cite{Ner-field}.

Oligomorphic groups $\GL(\infty)$ over finite fields also can be included
to a wider picture (see, e.g., \cite{Ner-finite}).

 So, there is a hope that  some other oligomorphic groups
 can be included to wider theories in representation-theoretical
 spirit.
 
 \sm
 
 4) Unitary representations of
 inverse limits $\lim\limits_\leftarrow G_n$ of oligomorphic groups 
 are easily reduced to representation of prelimit oligomorphic groups $G_n$, see \cite{Tsa}.
 Precisely, any irreducible representation skips through one of groups $G_n$
 (moreover, this is a general statement for any inverse limits of groups, this follows from \cite{Ner-book}, Proposition VIII.1.3).
  But here we get new possibilities. 
  
For example, consider the Bruhat--Tits tree $\cT_\infty$  of infinite order, 
i.e., the simplicial  tree whose vertices have countable valences.
Fix a vertex $\xi$ and consider a 'ball' $(\cT_\infty^n,\xi)\subset \cT_\infty$ of radius $n$
with center at $\xi$. Consider the group $\Aut(\cT_\infty)$ of all
automorphisms of $\cT_\infty$, its subgroup $\Aut(\cT_\infty,\xi)$ fixing $\xi$,
and the groups $\Aut(\cT_\infty^n,\xi)$. The groups $\Aut(\cT_\infty^n,\xi)$
are oligomorphic, the group $\Aut(\cT_\infty,\xi)$ is their
inverse limit,
$$
\Aut(\cT_\infty,\xi)=\lim_\leftarrow \Aut(\cT_\infty^n,\xi). 
$$

 Tsankov theorem \ref{th:tsankov} immediately
gives the classification of unitary representations of groups
$\Aut(\cT_\infty^n,\xi)$. 

An object of the corresponding category is
a finite rooted tree
$(L,\xi)$
such that distances between any vertex of $L$ and the root $\xi$ are $\le n$.
A morphism $(L,\xi)\to (M,\xi)$
is   an isomorphism of a rooted subtree $(A,\xi)\subset (L,\xi)$
to a rooted subtree $(B,\xi)\subset (M,\xi)$.

\sm 

We get a more rich situation, when we pass to the natural overgroup
of  $\Aut(\cT_\infty,\xi)$, namely, to $\Aut(\cT_\infty)$,
see%
\footnote{The train of $\Aut(\cT_\infty)$ is not the category of finite trees and their partial
isomorphisms as one might think looking to the previous examples, see \cite{Ner-book}, Sect VIII.6.} \cite{Olsh-new}, and to the group of spheromorphisms
of $\cT_\infty$, \cite{Ner-trees}.

\sm

Another examples are groups 
$\GL(\infty,\Z_p)$ and $\GL^\diamond(\infty,\Z_p)$
over the ring $\Z_p$ of $p$-adic integers, they are inverse limits
\begin{equation}
\GL(\infty,\Z_p)=\lim_\leftarrow \GL(\infty,\Z/p^n\Z),\quad
\GL^\diamond(\infty,\Z_p)=\lim_\leftarrow \GL^\diamond(\infty,\Z/p^n\Z).
\label{eq:p-adic}
\end{equation}
of oligomorphic groups. For the group $\GL(\infty,\Z/p^n\Z)$ objects of the train
are finite modules $L$ over the ring $\Z/p^n\Z$, a morphism $L\to M$ is an
isomorphism of a submodule $A\subset L$ to a submodule $B\subset M$.
For the  groups $\GL^\diamond(\infty,\Z/p^n\Z)$ explicit description of 
the trains is unknown, see \cite{Ner-p-adic}.   

Clearly, the groups \eqref{eq:p-adic} have natural overgroups
consisting of $p$-adic matrices of the same structure.

  \sm
  
%  4)
% In all our examples except  $G=S_\infty$
% double cosets  really are reduced. 

 \section{The group of automorphisms of the universal cubic spaces}
 
 \COUNTERS
  
  {\bf\punct The universal qubic space.} 
   Let $\F$ be  a finite  field with $q$ elements. We
   use the term {\it cubic space} for a linear space $V$  over $\F$ equipped with a symmetric trilinear ('cubic')
    form $\varpi$,
$$
\varpi: V\times V\times V \to \F.
$$      
We say that an invertible operator $V\to V$ is {\it isocubic} if it preserves
the cubic form. Denote $\Isoc(V)$ or $\Isoc(V,\varpi)$
 the group of all such maps.

     Let $V$, $W$ be qubic spaces equipped with
forms $\varpi_{\fan V}$, $\varpi_{\fan W}$.    
    We say that a linear injective
   map $\iota:V\to W$ is an {\it isocubic embedding} if for any $v_1$, $v_2$, $v_3\in V$
   we have 
   $$
   \varpi_{\fan W}\bigl(\iota(v_1), \iota(v_2), \iota(v_3)\bigr)=
   \varpi_{\fan V}(v_1,v_2,v_3).
   $$
   By $\St(V,W)$ (a {\it Stiefel variety}) we denote the set of all isocubic embeddings
   $V\to W$. The group $\Isoc(V)\times \Isoc(W)$ acts on $\St(V,W)$ in a natural way:
    $\iota \mapsto \sigma\iota \lambda^{-1}$, $\sigma\in \Isoc(W)$, 
   $\lambda\in \Isoc(V)$.
     
   We say that  a cubic space  $\cZ$ of countable dimension
  satisfies  the 
 {\it Urysohn extension property} if:
 
\sm 
 
 --- {\it  for any finite-dimensional cubic space
 $W$ and its subspace $V$ of codimension 1,
 any isocubic embedding $V\to \cZ$ extends to an isocubic embedding $W\to 
 \cZ$.}
 
\begin{theorem}
\label{th:universal}
{\rm a)} There exists a unique up to  linear equivalence
cubic space $\cZ$ of countable dimension 
satisfying the Urysohn extension property.

\sm

{\rm b)} The space $\cZ$ 
satisfies properties:

\sm

--- {\rm Universality}: any finite-dimensional cubic space admits an isocubic embedding
in $\cZ$;

\sm 

--- {\rm Ultrahomogenity}:  For any finite-dimensional
 cubic space $V$ the Stiefel variety $\St(V,\cZ)$ is $\Isoc(\cZ)$-homogeneous.
\end{theorem}

In particular, $\Isoc(\cZ)$ is oligomorphic, this follows from ultrahomogeneity.

\sm 

{\bf\punct Remarks.}
1) 
Theorem \ref{th:universal} is a very special case of Harman, Snowden \cite{HS}, on cubic spaces
see also \cite{KZ}, \cite{Bik1}, \cite{Bik2}. In our degree of generality, the statement admits a trivial probabilistic proof,  we need a probabilistic construction and present it in the next subsection.

\sm

2) Groups preserving symmetric or skew-symmetric bilinear forms 
are  subjects of substantial theories. Spaces with trilinear forms are a kind of 
exotics. However, exceptional Lie (and Chevalley) groups admit realizations as groups preserving (symmetric or skew-symmetric) trilinear forms, see, e.g., \cite{Gan}, \cite{Asch}, \cite{Wil}. 

\sm 

3)
Infinite-dimensional manifolds
of embeddings of a compact manifolds $M$ in a manifold $N$ can be equipped with
natural
  closed $k$-forms admitting infinite-dimensional groups of symmetries with infinite-dimensional stabilizers of points, see Vizman \cite{Viz}.

%\sm 

%4) Universal cubic spaces (and more generally universal tensor spaces)
%are representatives of numerous counterparts of the Urysohn universal metric
% spaces \cite{Ury} (as the Rado graph, the Hall universal group, the universal poset, the Gurarii space, ultrametric Urysohn spaces), see, e.g., \cite{Macph}.

\sm

{\bf \punct Existence of the universal cubic space.}
A. Consider a linear space $Z$ over $\F$ of countable dimension, fix a basis $e_1$, $e_2$, \dots. For each non-ordered triple of vectors we assign 
a random value
$\zeta(e_i,e_j,e_k)\in \F$ uniformly distributed in $\F$. We assume that all these values are independent.
 We get a random cubic form $\zeta$
on
 $V$. We also can say that we consider a countable direct product of copies of  the field $\F$
 equipped with a uniform measure, factors are enumerated by non-ordered triples of basis elements $(e_i,e_j,e_k)$, which can coincide.
 
 Equivalently, we  take an arbitrary  	exhaustive sequence
$X_1\subset X_2\subset\dots$ of finite dimensional subspaces such that
each $X_j$ is  equipped with a uniformly distributed cubic form.
 A restriction of  forms from $X_j$ to $X_{j-1}$ sends uniform distributions to
 uniform distributions; we take the inverse limit.
 
 Consider another exhaustive sequence $Y_1\subset Y_2\subset\dots$.
 Then each $Y_i$ is contained in some $X_{\alpha(i)}$ and each 
 $X_j$ is contained in some $Y_{\beta(j)}$. So it produces the same 
 distribution of cubic forms on $Z$.

  In particular, the distribution of forms does not depend on a choice of a basis.
 
\sm 
 
B. We claim that {\it with probability 1 a form $\zeta$ satisfies the Urysohn extension property}. Indeed, consider a finite dimensional cubic space $W$
 and its subspace $V$ of codimension 1. Fix an isocubic embedding 
 $\iota: V\to(Z,\zeta)$. Fix  a basis
 $e_1$, \dots, $e_k$ in $\iota(V)$ and complete it to a basis $\{e_\alpha\}$ in $Z$.
 For each $\alpha>k$ we consider the subspace $W_\alpha$ generated by $\iota(V)$
  and $e_\alpha$.
 The cubic form on $W_\alpha$ is determined by a collection of 
 $k(k+1)/2+k+1	$
  numbers 
 \begin{equation}
 \zeta(e_i, e_j,e_\alpha), \quad \zeta(e_i,e_\alpha,e_\alpha), \quad \zeta(e_\alpha,e_\alpha, e_\alpha)\in \F,
 \label{eq:collections}
 \end{equation}
  where $i$, $j\le k$. With probability 
 $q^{-(k+1)(k+2)/2}$ this space is isocubic to $W$. Since for different $\alpha$  collections
 \eqref{eq:collections} are independent, with probability $1$ we  meet a subspace
 isocubic to $W$. 
 
 The set of possible pairs of subspaces $V\subset W$ is countable, the set of  embeddings $V\to Z$ is countable,  so with probability 1 the extension of an embedding
 $V\to\cZ$ is possible for  any pair $V\subset W$. The Urysohn extension property
 is proved.
 
 \sm

C. Now the {\it universality} is obvious -- we take a basis $f_1$, \dots, $f_k$ in a cubic space $Y$, denote by $Y_j\subset Y$ the linear span of $f_1$, \dots, $f_j$. By the Urysohn property we can extend an
 isocubic embedding $Y_j\to \cZ$ to an isocubic embedding
  $Y_{j+1}\to \cZ$. We start with $Y_0=0$ and come to $Y_k=Y$.
  
\sm
  
D. Let us verify the {\it ultrahomogeneity}.  
Let $P_k$, $Q_k$ be two $k$-dimensional isocubic subspaces in a space 
$\cZ$ satisfying
the Urysohn property.
Consider a sequence $a_j$ containing once all elements of 
$\cZ\setminus P_k$
and a sequence $b_j$ containing once all elements of 
$\cZ\setminus Q_k$.

% Choose corresponding  bases $p_1$, \dots, $p_k\in P$
%and $q_1$, \dots, $q_k\in Q$. Complete them to bases $p_\alpha$ and $q_\beta$ in $\cZ$ in an arbitrary way.
 Denote $P_{k+1}:=P_k\dotplus \F a_{1}$.
 We extend the map $\iota_k:P_k\to Q_k$ to an isocubic embedding 
$\iota_{k+1}:P_{k+1}\to \cZ$. Denote the image by $Q_{k+1}$.
Remove from the sequence $a_j$ all elements that are contained 
in $P_{k+1}$ and shift numeration of remaining elements
preserving their order. Perform the same procedure with $Q_{k+1}$
and the sequence $b_k$.

% denote
%$\wt q_{k+1}=\iota_{k+1} p_{k+1}$.
%Vectors 
%$$
%q_1,\, q_2,\dots, q_k, \wt q_{k+1},\, q_{k+1}, q_{k+2},\dots
%$$
%are linear dependent. Remove the first element of this sequence
%that is a linear combination of the previous elements.
%Shifting numeration, we denote the remaining sequence by $q_1$, $q_2$, $q_3$,
%$\dots$ 
% Now new $q_{k+1}$, old $q_{k+1}$, and other $q_m$ become linear dependent. We remove the first $q_l$ such that $q_1$, \dots, $q_l$ are linear dependent and rename $q_s\mapsto q_{s-1}$ for $s>l$.
 Now we denote
$Q_{k+2}=Q_{k+1}\dotplus \F b_{1}$ (with new $b_1$),   take an isocubic embedding 
$\iota_{k+2}^{-1}:Q_{k+2}\to \cZ$, and denote the image by $P_{k+2}$. 
  We repeat the same operation a countable number of times, interlacing $P_m$ and $Q_m$, and come to an isocubic bijection $\cZ\to \cZ$ sending $P_k\to Q_k$.

\sm   
 
E. It remains to verify the {\it uniqueness} of $\cZ$. We take two cubic spaces $\cZ_1$, $\cZ_2$ satisfying the Urysohn property and repeat the  arguments from the proof of the ultrahomogeneity starting
from $P_0=0$ in $\cZ_1$ and $Q_0=0$ in $\cZ_2$.
\hfill $\square$ 
%
% assuming that
%$P_0\subset \cZ_1$, $Q_0\subset \cZ_2$ are zero-dimensional subspaces. 

\sm 

{\sc Remark.} This proof is a copy of the probabilistic proofs of existence of 
the  Rado graph, see Erd\H{o}s, P., R\'enyi \cite{ER}
(see, also, \cite{Cam-Rad}).
%or existence of the Urysohn space, see \cite{Ver}. 
The proof of universality is a copy
of Urysohn's considerations.
\hfill $\boxtimes$

\sm

{\bf\punct Some corollaries from the  probabilistic construction.}
Let $L$ be a linear space over $\F$ of countable dimension
with a fixed basis $e_l$ equipped with a cubic form $\lambda$.
Let $M$ be another linear space over $\F$ of countable dimension
with a fixed basis $f_m$. Denote be $\frZ(L)$ the space of all cubic 
forms $\zeta$ on%
\footnote{We denote by $S\dotplus T$  a direct sum of linear spaces,
they are not 'orthogonal' in any sense.} 
$$Z:=L\dotplus M$$
 extending the form $\lambda$ on $L$.
A point $(Z,\zeta)$ of $\frZ(L)$ is uniquely determined by a collection of values
$$
\zeta(e_l,e_{l'},f_k),\quad \zeta(e_l,f_{k},f_{k'}),\quad
\zeta(f_{k},f_{k'} ,f_{k''}), 
$$ 
where $l$, $l'$, $k$, $k'$, $k''\in\N$. We define a probabilistic measure
on $\frZ(L)$ assuming that all these random values are uniformly distributed 
and independent.

\begin{proposition}
\label{pr:refinement}
Almost all spaces $(Z,\zeta)\in\frZ(L)$ are isocubic to $\cZ$.
\end{proposition}

{\sc Proof.} It is sufficient to verify the Urysohn extension property.
Let $V\subset W$ be finite dimensional cubic spaces, $V$ has codimension 1 in W.
Consider an embedding $\phi:V\to L\dotplus M$, such that for any $v_1$, $v_2$, $v_3\in V$
satisfying $\phi(v_i)\in L$ we have
$$
\lambda\bigl(\phi(v_1),\phi(v_2),\phi(v_3)\bigr)=\varpi_V(v_1,v_2,v_3).
$$ 
Consider the subset $\frZ_\phi\subset\frZ$ consisting of all spaces $(Z,\zeta)$
such that 
$$
\zeta\Bigr|_L=\lambda, \qquad \zeta\bigl(\phi(v_1),\phi(v_2),\phi(v_3)\bigr)=
\varpi_V(v_1,v_2,v_3).
$$
Clearly, $\frZ_\phi$ is a subset of nonzero measure in $\frZ(L)$.
So, almost all spaces $(Z,\zeta)\in \frZ_\phi$ are isocubic to $\cZ$.

Consider a basis $h_\alpha\in V$, let $h_{\mathrm {next}}\in W\setminus V$.
Decompose vectors $\phi(h_\alpha)$ in the basis,
$$\phi(h_\alpha)=\sum a_{\alpha,i} e_i+\sum b_{\alpha,j} f_j.$$
Only finite number of $f_j$ actually enters to such sums, so all such $j$
are less than some $N$. We take a sequence of subspaces
$\phi(V)+\F f_N$, $\phi(V)+\F f_{N+1}$, $\dots$
All random values 
$$
\zeta\bigl(\phi(h_\alpha),\phi(h_\beta), f_{N+k}\bigr),
\quad \zeta\bigl(\phi(h_\alpha), f_{N+k},  f_{N+k}\bigr),
\quad \zeta\bigl(f_{N+k}, f_{N+k},  f_{N+k}\bigr)
$$
are independent. Therefore, a.s. there exists $k$ such that the linear map
$\wt\phi:W\to (Z,\zeta)$  defined  by
$\wt\phi(h_\alpha)=\phi(h_\alpha)$, $\wt \phi (h_{\mathrm {next}})=f_{N+k}$ is isocubic.
\hfill $\square$

\sm

We need the following modification of this proposition
(extracted from its proof).

Let $A\subset L$ be a subspace consisting of linear combinations 
of some basis vectors $e_{\alpha_1}$, \dots, $e_{\alpha_\mu}$. Let
$B\subset M$ be a subspace consisting of linear combinations 
of some basis vectors $f_{\beta_1}$, \dots, $f_{\beta_\nu}$.
Fix a cubic form $\varpi_{\vphantom{\bigl|}A\dotplus B}$
on $A\dotplus B$
 such that 
$\varpi_{\vphantom{\bigl|}A\dotplus B}\Bigr|_{A}=\lambda\Bigr|_A$.
Consider the space $\frZ\bigl(L|A,B\bigr)$ of all cubic forms $\zeta$
on $L\dotplus M$ such that 
$$
\zeta\Bigr|_{A\dotplus B}=\varpi_{\vphantom{\bigl|}A\dotplus B}, \qquad \zeta\Bigr|_L=\lambda.
$$ 
Such form is uniquely determined by collection
of values 
\begin{equation}
\zeta(e_i,e_{i'},f_j),\quad \zeta(e_i,f_j, f_{j'}),\quad \zeta(f_j, f_{j'}, f_{j''}),
\label{eq:random-values}
\end{equation}
where at least one $i$, $i'\notin \{\alpha_1,\dots,\alpha_\mu\}$
or at least one $j$, $j'$, $j''\notin \{\beta_1,\dots,\beta_\nu\}$.
We define the probabilistic measure on $\frZ\bigl(L|A,B\bigr)$
assuming that all random values
\eqref{eq:random-values} are independent and uniformly distributed.

\begin{corollary}
	\label{cor:refinement}
Almost all spaces $(Z,\zeta)\in\frZ\bigl(L|A,B\bigr)$ are isocubic to $\cZ$.
\end{corollary}

{\sc Proof.} The space $\frZ\bigl(L|A,B\bigr)$ is a subset of nonzero measure in 
$\frZ(L)$. \hfill $\square$

\sm

{\bf \punct Open subgroups of the group $\Isoc(\cZ)$.}
%For a space $W$ with cubic form denote by $\Isoc(W)$ the group of all isocubic
%bijections of $W$. 
For any finite-dimensional subspace $V\subset \cZ$ we denote by $K(V)\subset\Isoc(\cZ)$
the subgroup fixing all elements of $V$, by $K^\circ(V)\subset \Isoc(\cZ)$
the subgroup sending $V$ to itself. Clearly, 
$$
K^\circ(V)/K(V)\simeq \Isoc(V).
$$
The group $\Isoc(\cZ)$ is a closed subgroup in the symmetric group $S(\cZ)$,
the subgroups $K(V)\subset \Isoc(\cZ)$  are open
and form a fundamental system of neighborhoods of the unit.

\begin{proposition}
\label{pr:subspaces-perfect}
All finite-dimensional  subspaces in $\cZ$ are perfect subsets. 
\end{proposition}

\begin{theorem}
\label{pr:lemma}
For any open subgroup $H\subset\Isoc(\cZ)$ there is a
finite-dimensional subspace $V$ such that
$K(V)\subset H\subset K^\circ(V)$. 
\end{theorem}

Thus,  the subgroup $\Isoc(\cZ)\subset S(\cZ)$ is perfect oligomorphic, 
perfect subsets are linear subspaces in $\cZ$.
So, Theorem 
\ref{th:tsankov} gives us a classification of unitary representations
of $\Isoc(\cZ)$.
By Theorem \ref{th:main}, we have a well-defined category of reduced double cosets, it is equivalent to the category of partial linear isocubic bijections.

\begin{corollary}
\label{cor:category}
The category $\Tra\bigl(\Isoc(\cZ)\bigr)$
is equivalent to the following category: its objects 
 are finite-dimensional cubic linear spaces. A
morphism $V\to W$ is an isocubic linear bijection of a 
linear subspace $A\subset V$ to a linear subspace $B\subset W$.
\end{corollary}

\sm

{\bf\punct Proof of Proposition \ref{pr:subspaces-perfect}.}
%The statement is a corollary of the following lemma:

\begin{lemma}
\label{l:100}
Let $Y$, $V\subset\cZ $ be  finite-dimensional subspaces, $V\cap Y= 0$. 
Then
there exists a countable collection of
$g_j\in K(V)$ such that
the  sum 
$$L:=V+Y+ g_1 Y+g_2 Y+\dots$$
is direct and has countable codimension in $\cZ$.
\end{lemma}

 Proposition \ref{pr:subspaces-perfect} follows from the lemma. It is sufficient
 to consider one-di\-men\-si\-o\-nal $Y$.
 
\sm 

%By this lemma, any finite-dimensional subspace in $\cZ$ is a perfect subset
%in $\cZ$ 
%(Proposition \ref{pr:subspaces-perfect}).

%\sm 

{\sc Proof of Lemma \ref{l:100}.}
We apply Proposition \ref{pr:refinement}. Let $Y^1$, $Y^2$, $\dots$
be disjoint copies of $Y$.
Consider a direct sum 
\begin{equation}
L:=V\dotplus Y\dotplus Y^1\dotplus Y^2\dotplus\dots
.
\label{eq:YYY}
\end{equation}
We define cubic form on each subspace $V\dotplus Y\subset L$,
 $V\dotplus Y^j\subset L$
assuming that tautological maps from the initial space $V\dotplus Y$ 
to the corresponding subspace $V\dotplus Y\subset L$ and the subspaces $V\dotplus Y_j\subset L$
are isocubic. We consider an arbitrary cubic form on $L$
compatible with the forms on $V\dotplus Y$, $V\dotplus Y^j$.
Next, we consider the corresponding space $\frZ(L)$
and take an element $(Z,\zeta)\in \frZ(L)$ isomorphic to $\cZ$.
By the ultrahomogeneity, we have $g_j\in \Isoc(Z,\zeta)$ satisfying
desired conditions.
\hfill $\square$

\sm

{\bf\punct Proof of Theorem \ref{pr:lemma}.} Let $H$ be an open subgroup, $H\supset K(V)$, and there exists 
$\xi\in H\setminus K^\circ(V)$. It is sufficient to prove that there exists
a proper subspace $W\subset V$ such that $K(W)\subset H$.

Notice that $\xi K(V)\xi^{-1}=K(\xi V)$ and $\xi V\ne V$. 
So it is sufficient to prove the following statement
(assuming $V':=\xi V$):

\begin{lemma}
\label{l:1}
 Let a subgroup $H\subset \Isoc(\cZ)$ contain $K(V)$ and $K(V')$. Then it contains 
$K(V\cap V')$.
\end{lemma}

\begin{lemma}
\label{l:2}
In the notation of Lemma {\rm \ref{l:1}}, denote $W:=V\cap V'$.
 Decompose $V=W\dotplus X$, $V'=W\dotplus Y$.
Let $W\dotplus \wt X\dotplus \wt Y$ be a subspace in $\cZ$
such that there exists an isocubic map 
$r:W\dotplus X\dotplus Y\to W\dotplus \wt X\dotplus \wt Y$ 
identical on $W$ and
sending $X\to\wt X$, $Y\to \wt Y$. 
%Let 
%\begin{equation}
%(\wt X \oplus \wt Y) \cap (W\oplus X\oplus Y)=0
%\label{eq:for-omit}
%\end{equation}
 Then there is an element  $h\in H$ such that the restriction of $h$ to 
$W\dotplus X\dotplus Y$ is $r$.
\end{lemma}

{\sc Proof of Lemma   \ref{l:2}.}
We have $K(W\dotplus X)$, $K(W\dotplus Y)\subset H$.
Notice that
$$
h\in H,\quad K(A)\subset H\qquad \Rightarrow \qquad K(hA)=hK(A)h^{-1}\subset H.
$$

{\it Step 1.}
Let us show that {\it without loss of generality we can assume that} 
\begin{equation}
\wt X\cap \bigl(W\dotplus X\dotplus Y\bigr)=0.
\label{eq:transversal}
\end{equation}

By Lemma \ref{l:100} and ulrahomogeneity, there is a countable collection of maps 
$g_j\in K(W\dotplus Y)$ such that a sum 
$$W+Y+X+g_1 X+ g_2 X+\dots$$
 is direct. For sufficiently large $N$, we have
$\wt X\cap (W\dotplus g_N X)=0$.
Now we can replace $X\mapsto g_N X$, $r\mapsto r g_N^{-1}$.
So, we can assume 
$$\wt X\cap (W\dotplus X)=0.$$

We again apply Lemma \ref{l:100}. There exists a countable collection of maps
$h_j\in K(W\dotplus X)$ such that  a sum
$W+X+Y+h_1 Y+ h_2 Y+\dots$ is direct. For sufficiently large $N$, we 
have $\wt X\cap (W\dotplus X\dotplus h_N Y)=0$.
Now we can replace $Y\mapsto h_N Y$, $r\mapsto r h_N^{-1}$.
We come to the case
\eqref{eq:transversal}.
 
\sm

{\it Step 2.} So, we assume that \eqref{eq:transversal} holds.
 Now we intend to apply
   Lemma \ref{l:100} and Corollary \ref{cor:refinement} again.
Choose bases 
\begin{equation}
w_\alpha\in W,\quad x_\beta\in X, \quad y_\gamma\in Y.
\label{eq:choose-bases}
\end{equation}
Consider  the space $L$ as in the proof of Lemma \ref{l:100} assuming $V=W\dotplus X$. Consider copies $y_\gamma^j\in Y^j$
 of vectors $y_\gamma$. Such vectors together with \eqref{eq:choose-bases} form a basis in the space $L$.
 So we get a space $L$ equipped with a cubic form, say $\lambda$.
 
 Next, consider the space $Z:=L\dotplus M$ as in Proposition \ref{pr:refinement}. Let us identify
 the basis $rx_\alpha\in\wt X$ with first basis elements $f_1$,  \dots, $f_s \in M$.
 Consider the restriction of our initial cubic form $\zeta$ to $W\dotplus X\dotplus Y\dotplus \wt X$,
 and assign this cubic form to the corresponding subspace  in $L\dotplus M$.
 
 We came to an object described in Corollary \ref{cor:refinement},
 we have fixed cubic forms on $L$ and $W\dotplus X\dotplus Y\dotplus \wt X$
 and these forms are compatible on their intersection  $W\dotplus X\dotplus Y$.
 So, we get a probabilistic space $\frZ(L|W\dotplus X\dotplus Y, \wt X)$
 whose points are forms $\zeta$ on $Z:=L\dotplus M$ extending our forms to the whole
space. By Corollary \ref{cor:refinement}, a cubic space $(Z,\zeta)$ a.s. is isomorphic
to $\cZ$. In particular, this space contains a copy of $\wt Y$ such that
a new $W\dotplus X \dotplus Y\dotplus \wt X+ \wt Y$ is isocubic to the initial
$W\dotplus X \dotplus Y\dotplus \wt X+ \wt Y$. Position  of new $\wt Y$
depends on  a random form $\zeta$,
 choice of this $\wt Y$ is not unique (it can be done in a canonical way, but we do not need this).

% Let us extend each $p_j$ to an element 
% $g_j\in K(V)=K(W\oplus X)$.
%Choose a collection $g_j$ as in Lemma \ref{l:100} taking $V=W\dotplus X$.
% Notice that $K(W\oplus g_j Y)=g_j K(W\oplus Y)g_j^{-1}\subset H$.  

%Denote by $Q_{\vphantom{\bigl|}\cZ}$ the trilinear form on $\cZ$. 
 Let
 $j$ range in  $1$, $2$, $3$,  \dots, and $w_\alpha$, $x_\beta$, $y_\gamma$
 range in sets of basis vectors of $W$, $X$, $Y$ respectively.
  Then the following collection of $\F$-valued random variables
$$
\zeta( y_\gamma^j,  y_{\gamma'}^j, r x_{\beta}),\quad
\zeta( y_\gamma^j, r x_{\beta}, r x_{\beta'}),\quad
\zeta( y_\gamma^j, w_\beta, r x_{\beta})
$$
is independent. Therefore, with probability 1 there is $j_*$ such that
the following (finite) collection of identities holds:
\begin{align}
\zeta(y_\gamma^{j_*},  y_{\gamma'}^{j_*}, r x_{\beta})&=
\zeta( y_\gamma^{j_*}, y_{\gamma'}^{j_*},  x_{\beta});
\label{eq:1}
\\
\zeta( y_\gamma^{j_*},  r x_{\beta}, r x_{\beta'})&=
\zeta( y_\gamma^{j_*},   x_{\beta},  x_{\beta'});
\label{eq:2}
\\
\zeta( y_\gamma^{j_*}, w_\beta, r x_{\beta})&=
\zeta( y_\gamma^{j_*}, w_\beta, 
x _{\beta}).
\label{eq:3}
\end{align}
Consider the map 
$$\sigma:W\dotplus X\dotplus Y^{j_*}\to W\dotplus rX\dotplus  Y^{j_*},$$
which is identical on $W\dotplus Y^{j_*}$ and equals
 $r$ on $X$ (and therefore $\sigma$ equals $r$ on $W\dotplus X$). This map is  
 $\zeta$-isocubic. Indeed,
 
 \sm
 
 1) by  definition, it is isocubic
 on $W\dotplus  Y^{j_*}$;
 
\sm 
 
 2) $r$ is isocubic on $W\dotplus X$;

\sm 
 
 3) by \eqref{eq:1}-\eqref{eq:2}
 it is isocubic on $X\dotplus  Y^{j_*}$;
 
\sm 
 
 4) by \eqref{eq:3} it preserves
 the trilinear form on $ Y^{j_*}\times X\times W$.
 
\sm 
 
   The map $\sigma$ extends to
 an element $\wt\sigma$ of the subgroup $K(W\dotplus  Y^{j_*})\subset H$. 
Conjugating $K(W\dotplus X)$ by $\wt\sigma$, we get that 
$K(W\dotplus \wt X)\subset  H$. 
Also, by the construction,
we have the following isocubic maps
$$
W\dotplus \wt X\dotplus \wt Y
\,\leftarrow\,
W\dotplus  X\dotplus  Y
\,\rightarrow\,
W\dotplus  X\dotplus  Y_{j^*} 
\,\rightarrow\,
W\dotplus  \wt X\dotplus  Y_{j^*}.
$$ 
So, we have an isocubic map
$W\dotplus \wt X \dotplus  Y^{j_*}\to W\dotplus \wt X \dotplus \wt Y$.
By  the construction, it fixes $W\dotplus \wt X$.
Therefore there is an element $\xi$ of $K(W\dotplus \wt X)$ sending one subspace to another. We take $h=\xi \wt\sigma$ \hfill $\square$

\sm

{\sc Proof of Lemma \ref{l:1}.}
Let $r\in H$. By Lemma \ref{l:2}  there exists $h\in H$
such that $h=r$ on $V+V'=W\dotplus X\dotplus Y$. So 
$h^{-1}r\in K(V+V')$. But $K(V+V')\subset K(V)$.
\hfill
$\square$

\sm

{\bf \punct Compactification of $\Isoc(\cZ)$.}
Consider the following semigroup $\Gamma(\cZ)$. 
Its elements are partial linear isocubic bijections
from a subspace $A\subset \cZ$ to a subspace $B\subset \cZ$.
The product is the product of partial bijections.
A sequence $\gamma_j\in \Gamma(\cZ)$ converges to $\gamma\in \Gamma(\cZ)$
if for each $v$, $w\in\cZ$ the condition '$\gamma v=w$'
is equivalent to the condition '$\gamma_j v=w$
for sufficiently large $j$'.

\begin{corollary}
{\rm a)}
The  $\Gamma(\cZ)$ is a   compact semigroup with separately continuous
product.

\sm 

{\rm b)} The group $\Isoc(\cZ)$ is dense in $\Gamma(\cZ)$.

\sm 

{\rm c)} Each unitary representation 
of $\Isoc(\cZ)$ extends to a weakly continuous continuous representation 
of $\Gamma(\cZ)$. 
\end{corollary}

In particular, for each unitary representation $\rho$ of 
the group
$\Isoc(\cZ)$, the closure $\ov{\rho\bigl(\Isoc(\cZ)\bigr)}$ of
$\rho\bigl(\Isoc(\cZ)\bigr)$ in the weak operator topology
coincides with the image of $\Gamma(\cZ)$.

\sm 

{\sc Proof.} a) The semigroup of all partial bijections 
of a countable set is
a compact semigroup with separately continuous multiplication,
see, e.g., \cite{Ner-book}, Sect. VIII.1. It is easy
to see that $\Gamma(\cZ)$ is closed in the semigroup of all 
partial bijection.

\sm 

b) Consider a partial isocubic bijection $\theta\in\Gamma(\cZ)$,
let $V:=\dom \theta$. Consider a sequence of  
finite-dimensional subspaces
$V_1\subset V_2\subset V_3\subset \dots $ such that $V=\cup V_j$.
Consider a sequence $x_1$, $x_2$, \dots consisting of all elements
of the set $\cZ\setminus V$. We choose a sequence $g_k\in \Isoc(\cZ)$
 such that for each $k$:
 
 \sm
 
 --- $g_k\Bigr|_{V_k}= \theta\Bigr|_{V_k}$;
 
\sm 
 
---  for each $j<k$, $i<k$ we have $g_k x_i\ne g_j x_i$.

\sm

This is possible due to Neumann Lemma \ref{l:finite-shift}. So, for each $x_m\in \cZ\setminus V$
the sequence $g_j x_m$ consists of pairwise different elements.
On the other hand, each $v\in V$ is contained in $V_k$ with sufficiently
large $k$, so $g_k v=\theta v$. Therefore, $g_k$ converges to
$\theta$.

\sm 

c) It is sufficient to verify this for irreducible representations.
We use the model of irreducible representations as Subsect. \ref{ss:rep-tra}.
Let $A\subset \cZ$ be a finite-dimensional subspace, $\sigma$
an irreducible representation of $\Isoc(A)$. The space 
$H(A,\sigma)$ is defined as Subsect. \ref{ss:rep-tra}. 
Operators $\tau(\theta)$
 in the space $H(A,\sigma)$ are defined 
 by
 $$
\tau(\theta) F(\iota)=\begin{cases}
F(\theta^*\iota),\quad &\text{if}\quad A\subset \im \theta,\, 
\theta^*\iota A \subset \dom \theta;
\\
0,\quad &\text{otherwise.}
\end{cases}
$$

\sm

{\bf \punct Final remarks.}
 In our considerations, symmetric trilinear forms can be replaced by 
skew-symmetric trilinear form, and more generally by symmetric or skew-symmetric
$n$-linear forms. We also can equip a linear space with a finite collection
of  such forms.

 Apparently, our Theorem \ref{pr:lemma} can be extended to arbitrary
 contravariant irreducible tensors over finite fields as in \cite{HS} (a non-obvious place is an extension of Proposition \ref{pr:lemma}). Apparently, construction
of \cite{Ner-book}, Sect.III.3 allows to consider also exterior forms of middle dimension.

\tt 

University of Graz,
\\
\phantom{.}
\hfill Department of Mathematics and Scientific computing;

High School of Modern Mathematics MIPT,
\\
\phantom{.}
\hfill 1 Klimentovskiy per., Moscow; 

Moscow State University, MechMath. Dept;

 University of Vienna, Faculty of Mathematics.
 
 \sm

e-mail:yurii.neretin(dog)univie.ac.at

URL: https://www.mat.univie.ac.at/$\sim$neretin/

\phantom{URL:} https://imsc.uni-graz.at/neretin/index.html

\end{document}